\documentclass[11pt,a4paper]{amsart}
\usepackage[T1]{fontenc}
\usepackage{lmodern}
\usepackage{amssymb,amsmath,mathtools}
\usepackage[margin=2.6cm]{geometry}
\usepackage{xcolor}
\usepackage{tikz}
\usepackage{booktabs}
\usepackage[round]{natbib}
\usepackage{xspace}
\usetikzlibrary{calc}
\usepackage[colorlinks=true,linkcolor=blue!70!black,citecolor=blue!70!black,urlcolor=blue!70!black]{hyperref}

\newcommand{\doi}[1]{doi: \href{https://doi.org/#1}{\nolinkurl{#1}}}

\newtheorem{theorem}{Theorem}[section]
\newtheorem{lemma}[theorem]{Lemma}
\newtheorem{proposition}[theorem]{Proposition}
\newtheorem{corollary}[theorem]{Corollary}
\theoremstyle{definition}
\newtheorem{definition}[theorem]{Definition}
\newtheorem*{definition*}{Definition}
\newtheorem{example}[theorem]{Example}
\newtheorem{problem}{Problem}
\theoremstyle{remark}
\newtheorem{remark}[theorem]{Remark}

\newcommand{\Z}{\mathbb{Z}}
\newcommand{\F}{\mathbb{F}}
\newcommand{\Q}{\mathbb{Q}}
\newcommand{\C}{\mathbb{C}}
\newcommand{\Cay}{\operatorname{Cay}}
\newcommand{\EK}{\operatorname{EK}}
\newcommand{\wdd}{\operatorname{wd}_{2,1}}
\newcommand{\rab}{\rho_{\mathrm{ab}}}
\newcommand{\eall}{\eta_{\mathrm{all}}}
\newcommand{\ind}{\mathbf{1}}
\newcommand{\jv}[1]{\mathsf{#1}}   

\hypersetup{
  pdftitle={Degree deficiency, rigidity and a spectral obstruction for induced embeddings into abelian Cayley graphs},
  pdfauthor={Rigobert Fokam Souop and Laurent Bitjoka}
}

\makeatletter
\g@addto@macro\@setauthors{%
  \begin{center}\footnotesize
  \textit{Department of Electrical Engineering, Energy and Automation, ENSAI,}\\
  \textit{University of Ngaound\'er\'e, P.O.~Box 455, Ngaound\'er\'e, Cameroon}\\[4pt]
  \textit{E-mail addresses:} \href{mailto:fokamrigobert@gmail.com}{fokamrigobert@gmail.com} (R.~Fokam Souop),\\
  \href{mailto:laurent.bitjoka@univ-garoua.cm}{laurent.bitjoka@univ-garoua.cm} (L.~Bitjoka)
  \end{center}}
\makeatother

\begin{document}

\title[Degree deficiency, rigidity and a spectral obstruction]{Degree deficiency, rigidity and a spectral obstruction for induced embeddings into abelian Cayley graphs}

\author{Rigobert Fokam Souop}
\author{Laurent Bitjoka}


\begin{abstract}
Let $\eta(G)$ be the least order of a finite abelian group whose Cayley graph contains the connected graph $G$ as an induced subgraph. Generator classes give a weighted edge decomposition, whose minimum cost we characterize for arbitrary graphs. Commutativity gives stronger restrictions. If $G$ contains an induced copy of $K_4$ with one edge subdivided, then no host has degree four once $G$ has eight vertices, and none has degree five once $\eta(G)\ge15$; the second result rests on a finite computation. We also prove $\eta(G)\ge 2\lceil 7n/10\rceil$ for cubic graphs of girth at least five on $n$ vertices. The Erd\H{o}s--Kelly conditions and host-degree lower bounds alone give at most $14$, $28$ and $17$ for the Petersen graph, the Desargues graph and a specified eleven-vertex graph, respectively. A character-sum obstruction for tight regular embeddings gives a computer-free proof of $\eta=16$ for Petersen and $\eta\ge30$ for Desargues. SAT computations with checked proofs give $\eta=32$ for Desargues, $20$ for the eleven-vertex graph and $36$ for the smallest cubic graph without a perfect matching.
\end{abstract}

\thanks{Corresponding author: Rigobert Fokam Souop.}
\maketitle

\section{Introduction}

For a finite graph $G$, let $\eta(G)$ be the least order of a finite abelian group $\Gamma$ such that some Cayley graph $\Cay(\Gamma,S)$ contains $G$ as an induced subgraph. From Section~2 onward, $G$ is connected unless stated otherwise. Every graph has such a host \citep{Ba78,BS85}, and $\eta(G)=O(n^2)$ for graphs on $n$ vertices. In \citet{FB26c} we determined $\eta$ for complete multipartite graphs and for blow-ups, studied double stars and trees, and computed $\eta$ for all connected graphs on at most seven vertices. The lower bounds there come mainly from the conditions of \citet{EK63}: if a graph $G$ with $n$ vertices, $m$ edges and minimum degree $\delta$ is an induced subgraph of a $d$-regular graph on $n+t$ vertices, then
\begin{equation}\label{eq:EKintro}
td\ge nd-2m,\qquad t\ge d-\delta,\qquad t^2-(d+1)t+nd-2m\ge 0,\qquad (n+t)d\ \text{even}.
\end{equation}
A Cayley graph $\Cay(\Gamma,S)$ is $|S|$-regular, so \eqref{eq:EKintro} holds with $d=|S|$. The first two conditions are increasing in $d$, and in general the only information on the host degree is $d\ge\Delta(G)$, improved to $d\ge\Delta(G)+1$ for regular graphs that are not abelian Cayley graphs \citep[Lemma~2.5(v)]{FB26c}. This paper gives lower bounds on the host degree that use commutativity, shows where bounds of this kind stop, and goes past that point with a spectral argument.

\smallskip\noindent\textbf{Generator classes.} A connection set $S$ is a union of classes $\{s,-s\}$. If $2s\ne0$, the edges of the host with difference $\pm s$ form a disjoint union of cycles; if $2s=0$, they form a perfect matching. Restricted to an induced copy of $G$, the classes therefore partition $E(G)$ into parts of maximum degree at most two, each costing two elements of $S$, and matchings, each costing one. Hence $|S|\ge\wdd(G)$, where $\wdd(G)$ is the least value of $2a+b$ over all partitions of $E(G)$ into $a$ parts of maximum degree at most two and $b$ matchings. Colouring the edges of a Cayley graph by generator classes is standard; it underlies the no-lonely-colour colourings of minimal Cayley graphs studied by \citet{GMK25} and \citet{KSG25}.

\smallskip\noindent\textbf{Results.}
\begin{enumerate}
\item \emph{Rigidity of the subdivided $K_4$} (Theorems~\ref{thm:deg4} and~\ref{thm:deg5}). Let $J$ be $K_4$ with one edge subdivided. Every $4$-regular abelian Cayley graph containing an induced copy of $J$ has all its components isomorphic to $\Cay(\Z_7,\{\pm1,\pm2\})$. For degree five, a computer-assisted classification shows that the components have $8$, $10$ or $14$ vertices. Hence a connected graph $G$ containing an induced $J$ has no host of degree four if it has at least eight vertices, and no host of degree at most five if $\eta(G)\ge15$ (Corollary~\ref{cor:deg6}). At degree six the restriction disappears: $C_7^2\,\square\,C_m$ contains an induced $J$ for every $m\ge3$.
\item \emph{Commuting squares} (Proposition~\ref{prop:square} and Corollaries~\ref{cor:C4free} and~\ref{cor:girth5}). A $k$-regular graph without four-cycles needs abelian host degree at least $2k-1$ for odd $k$ and $2k-2$ for even $k$. For cubic graphs of girth at least five this gives $\eta(G)\ge2\lceil7n/10\rceil$; the degree bound is attained by the Desargues graph.
\item \emph{Generator classes in arbitrary groups} (Theorems~\ref{thm:wd} and~\ref{thm:nonab}). Let $X_\Delta$ be the set of vertices of maximum degree $\Delta$. Then $\wdd(G)=\Delta+1$ if $\Delta$ is odd and no matching of $G$ covers $X_\Delta$, and $\wdd(G)=\Delta$ otherwise; this is computable in polynomial time, whereas the chromatic index is NP-hard to compute already for cubic graphs \citep{Ho81}. The bound $|S|\ge\wdd(G)$ holds for Cayley graphs of all finite groups, so \eqref{eq:EKintro} holds with $d\ge\wdd(G)$ for every host, abelian or not. For regular graphs this adds nothing to \citet[Lemma~2.5(v)]{FB26c} (Remark~\ref{rem:regular}); for irregular graphs it can raise the Erd\H{o}s--Kelly floor (Example~\ref{ex:Bprime}).
\item \emph{Two extremal graphs} (Theorem~\ref{thm:2D5} and Proposition~\ref{prop:wallis}). The smallest connected graphs with $\wdd=\Delta+1$ have at least $2\Delta+5$ vertices. We prove the bound sharp for every odd $\Delta\ge3$; on eleven vertices there is exactly one example, the graph $B$ of Definition~\ref{def:B}; every host of $B$ has degree at least six. The smallest connected cubic graph without a perfect matching has sixteen vertices and is unique \citep{Wa81,CRW86}; for this graph $H$ and an infinite family of cubic graphs $G_T$ that includes it, every host has degree at least six, and $\eta(G_T)\ge3n/2$ (Theorem~\ref{thm:GT}).
\item \emph{The limit of degree arguments} (Proposition~\ref{prop:limit}). An argument that only bounds the host degree from below gives at most $\EK(G;\rab(G))$, where $\rab(G)$ is the least degree of a host. For the Petersen graph $P$, the graph $B$ and the Desargues graph $D$ our bounds $14$, $17$ and $28$ reach this limit; for $H$ no such argument can give more than $26$.
\item \emph{A spectral obstruction} (Theorems~\ref{thm:petersen} and~\ref{thm:desargues}). If a $k$-regular graph meets condition $td\ge nd-2m$ of \eqref{eq:EKintro} with equality, which we call a tight embedding (Section~\ref{sec:spectral}), the host vertices outside the image are independent and form one part of an equitable partition, so $k-d$ is an eigenvalue of the host and the indicator of the complement lives on the corresponding characters (Lemma~\ref{lem:tight}). With a Galois average of character sums (Lemma~\ref{lem:galois}) this classifies tight embeddings in groups of odd order with host degree $2k$: they are bipartite Haar graphs on two cosets of an index-three subgroup (Theorem~\ref{thm:odd}). It proves $\eta(P)\ge16$ by hand and $\eta(D)\ge30$, beyond the limit of the previous item.
\item \emph{Completions} (Theorem~\ref{thm:transfer}). A Cayley completion of $G$, a host of order $n$ containing $G$ as a spanning subgraph, has a generator partition too. A secondary labelling that separates, in every class, the edges of $G$ from the pairs the completion pays for turns it into an induced embedding of product form; for a binary labelling this is a linear system over $\F_2$.
\item \emph{Exact values} (Theorem~\ref{thm:exact}). By SAT solving, with a checked proof for every unsatisfiable instance, $\eta(P)=16$, $\eta(B)=20$, $\eta(D)=32$ and $\eta(H)=36$. The value $\eta(P)=16$ agrees with \citet[Table~5]{FB26c}.
\end{enumerate}

The conditions \eqref{eq:EKintro} are those of \citet{EK63}; Lemma~\ref{lem:petersen} is Petersen's $2$-factor theorem; and the order and uniqueness of the smallest cubic graph without a perfect matching are classical \citep{Wa81,CRW86}. We include a proof of the latter only because Theorem~\ref{thm:2D5} reuses its argument. Theorem~\ref{thm:wd} follows quickly from Petersen's theorem and a matching argument, and closely related objects have been studied. Such partitions are related to the generalized edge colourings of \citet{HK86}. \citet{Da26} considers cycle-matching colourings, in which every colour class is a disjoint union of isolated edges and odd cycles, and shows, using a theorem of \citet{CN13}, that $\Delta(G)$ such classes suffice for every multigraph. These colourings are unweighted: each class counts once, whatever its type. In $\wdd$ a part of maximum degree two costs two and a matching costs one, so for odd $\Delta$, the cost exceeds $\Delta$ exactly when no matching covers $X_\Delta$; we have not found this weighted version in the literature. The sharpness statement of Theorem~\ref{thm:sharp} is an additional general construction. Our other results are Theorems~\ref{thm:wd}, \ref{thm:deg4}, \ref{thm:deg5}, \ref{thm:2D5}, \ref{thm:GT}, \ref{thm:exact}, \ref{thm:odd}, \ref{thm:petersen}, \ref{thm:desargues} and~\ref{thm:transfer}, Lemma~\ref{lem:tight}, and Corollaries~\ref{cor:deg6}, \ref{cor:C4free} and~\ref{cor:girth5}. Proposition~\ref{prop:limit} is a simple observation, stated because it shows where degree arguments stop, and Lemma~\ref{lem:galois} is a standard Galois average, which we include for completeness.

\begin{remark}[Companion-paper scope]\label{rem:scope}
The companion study \citep{FB26c} concerns extremal host orders, trees, multipartite graphs, and exhaustive small-graph tables. The present article instead develops generator-class restrictions, a matching-saturation characterization, and abelian rigidity. The degree-five classification and the sharpness construction for every odd $\Delta$ do not follow from the companion paper. Computational infrastructure used for host-order searches is shared, but the structural proofs are independent.
\end{remark}

\smallskip\noindent\textbf{Related work.} The question of which graphs occur as induced subgraphs of Cayley graphs goes back to \citet{Ba78} and \citet{BS85}; \citet{Sp83} showed that there are graphs of arbitrarily high girth and bounded degree that are not induced subgraphs of any minimal (irreducible) Cayley graph. Representations modulo an integer \citep{EE89,EFMMP94} prescribe a cyclic group and a fixed set of allowed differences; \citet{GI87} study the existence of embeddings into Cayley graphs. \citet{KSG25} and \citet{GMK25} study subgraphs of minimal Cayley graphs, whose connection sets are inclusion-minimal generating sets; \citet{KSG25} show that some cubic graph of girth five is not a subgraph of any minimal Cayley graph, and that the generalized Petersen graphs $G(n,k)$ with $\gcd(n,k)=1$ are induced subgraphs of minimal Cayley graphs. Neither question concerns the least order of a host. The abelian Cayley result of \citet{PT20} concerns the maximum degree of induced subgraphs on more than half of the vertices. In the Cayley completion problem \citep{FB26a} the group has order $n$, every edge of $G$ is kept, and additional edges are paid for; there, as here, the generator classes are the basic structure, and Section~\ref{sec:completions} makes the connection precise.

\section{Preliminaries}

Graphs are finite and simple. Throughout, $G=(V,E)$ is connected with $n\ge2$ vertices and $m$ edges, unless a statement says otherwise. We write $d_G(v)$, $N_G(v)$, $\Delta=\Delta(G)$, $\delta=\delta(G)$, $\alpha(G)$ and $\chi'(G)$ for degrees, neighbourhoods, maximum and minimum degree, independence number and chromatic index. For $F\subseteq E$, $G[F]$ is the spanning subgraph $(V,F)$, so isolated vertices are kept. Except in Theorem~\ref{thm:nonab}, groups are finite and abelian and are written additively.

\subsection{Induced embeddings}
A \emph{connection set} of $\Gamma$ is a set $S\subseteq\Gamma\setminus\{0\}$ with $S=-S$. The Cayley graph $\Cay(\Gamma,S)$ has vertex set $\Gamma$, and distinct $x,y$ are adjacent when $x-y\in S$; it is $|S|$-regular. An element $s\ne0$ with $2s=0$ is an \emph{involution}.

\begin{definition}\label{def:embedding}
An \emph{induced embedding} of $G$ is a triple $(\Gamma,S,\lambda)$ of a finite abelian group $\Gamma$, a connection set $S$ and an injection $\lambda\colon V\to\Gamma$ such that, for distinct $u,v\in V$,
\[
uv\in E\iff \lambda(u)-\lambda(v)\in S.
\]
Its \emph{order} is $|\Gamma|$ and its \emph{degree} is $|S|$. The least order is $\eta(G)$.
\end{definition}

An induced embedding preserves adjacency and non-adjacency, but not necessarily distances; isometric embeddings, which also preserve distances, are not considered here.

For an injection $\lambda\colon V\to\Gamma$ put
\[
D_E(\lambda)=\{\lambda(u)-\lambda(v):uv\in E\},\qquad D_N(\lambda)=\{\lambda(u)-\lambda(v):u\ne v,\ uv\notin E\},
\]
where both orientations of each pair are included.

\begin{lemma}[Difference criterion]\label{lem:diff}
An injection $\lambda\colon V\to\Gamma$ is part of an induced embedding if and only if $D_E(\lambda)\cap D_N(\lambda)=\emptyset$. In that case $S=D_E(\lambda)$ is the smallest admissible connection set.
\end{lemma}

\begin{proof}
For an induced embedding, $D_E(\lambda)\subseteq S$ and $D_N(\lambda)\cap S=\emptyset$. Conversely, $D_E(\lambda)$ is symmetric and avoids $0$ because $\lambda$ is injective, so it is a connection set, and it separates edges from non-edges when $D_E\cap D_N=\emptyset$. Every admissible $S$ contains $D_E(\lambda)$.
\end{proof}

This is the abelian case of \citet[Proposition~3.1]{BS85}; see also \citet[Lemma~2.1]{FB26c}. Labelling $V$ by part of a Sidon set shows that $\eta(G)\le p^2$ for every odd prime $p\ge n$ \citep[Section~4]{BS85,FB26c}. Moreover $\eta(G)\ge n$, with equality exactly when $G$ is an abelian Cayley graph \citep[Lemma~2.3]{FB26c}.

\begin{lemma}[Connected reduction]\label{lem:connected}
Let $(\Gamma,S,\lambda)$ be an induced embedding of the connected graph $G$, and let $K=\langle S\rangle$. The components of $\Cay(\Gamma,S)$ are the cosets of $K$, each isomorphic to $\Cay(K,S)$, and $\lambda(V)$ lies in one of them. Translating gives an induced embedding into $\Cay(K,S)$; in particular some embedding of order $\eta(G)$ has $\langle S\rangle=\Gamma$.
\end{lemma}

\begin{proof}
A walk from $x$ changes the group element by a sum of elements of $S$, so it stays in $x+K$; conversely every element of $x+K$ is reached from $x$ by such a walk. Translation $y\mapsto y-x$ maps $x+K$ onto $K$ and preserves differences. Since $G$ is connected, any two labels are joined by the image of a path of $G$, so $\lambda(V)$ lies in one component, and translating it into $K$ preserves adjacency and non-adjacency.
\end{proof}

\subsection{The generator partition}

\begin{definition}\label{def:genpart}
Let $(\Gamma,S,\lambda)$ be an induced embedding. Let $\mathcal C(S)$ be the set of classes $\{s,-s\}$ with $s\in S$, and for $c\in\mathcal C(S)$ let $E_c=\{uv\in E:\lambda(u)-\lambda(v)\in c\}$. The family $(E_c)_{c\in\mathcal C(S)}$ is the \emph{generator partition}; it partitions $E$, with some parts possibly empty, and all parts are non-empty when $S=D_E(\lambda)$.
\end{definition}

\begin{lemma}\label{lem:classes}
Let $\lambda\colon V\to\Gamma$ be injective with $D_E(\lambda)\subseteq S$, and let $c\in\mathcal C(S)$.
\begin{enumerate}
\item[(a)] If $c=\{s\}$ with $2s=0$, then $G[E_c]$ is a matching together with isolated vertices.
\item[(b)] If $c=\{s,-s\}$ with $2s\ne0$, then $\Delta(G[E_c])\le2$, and every cycle of $G[E_c]$ has length $\operatorname{ord}(s)$.
\end{enumerate}
If $S$ has $a$ classes of size two and $b$ of size one, then $|S|=2a+b$.
\end{lemma}

\begin{proof}
Let $X_s$ be the spanning subgraph of the host with edges $\{x,x+s\}$, $x\in\Gamma$. If $s=-s$, then $x\mapsto x+s$ is a fixed-point-free involution and $X_s$ is a perfect matching; its restriction to $\lambda(V)$ is a matching, which proves (a). If $s\ne-s$, translation by $s$ permutes each coset of $\langle s\rangle$ cyclically, with orbits of length $r=\operatorname{ord}(s)\ge3$, so $X_s$ is a disjoint union of cycles $C_r$. The graph $G[E_c]$ is carried injectively into $X_s$, so its maximum degree is at most two, and a cycle in a subgraph of a disjoint union of cycles $C_r$ is one of those cycles. This proves (b). The last assertion holds because the classes partition $S$.
\end{proof}

The conditions of Lemma~\ref{lem:classes} are necessary but far from sufficient for a partition of $E$ to be a generator partition: the edge increments must also satisfy every cycle relation of $G$, the resulting labels must be distinct, and no non-edge difference may lie in $S$. Section~\ref{sec:wd} uses only the degree conditions; Section~\ref{sec:abelian} uses the relations between generators.

\subsection{Degree deficiency}

\begin{proposition}\label{prop:EK}
Let $G$ be an induced subgraph of a $d$-regular graph $X$ on $N$ vertices. Put $t=N-n$ and $\sigma_d=nd-2m=\sum_{v\in V}(d-d_G(v))$, let $U$ be the vertex set of the copy of $G$ and $W$ its complement. Then $X$ has exactly $\sigma_d$ edges between $U$ and $W$, and $2|E(X[W])|=td-\sigma_d$. Consequently
\begin{align}
\Delta\le d&\le N-1,\label{eq:c2}\\
td&\ge\sigma_d,\label{eq:c3}\\
t&\ge d-\delta,\label{eq:c4}\\
t^2-(d+1)t+\sigma_d&\ge0,\label{eq:c5}\\
Nd&\ \text{is even.}\label{eq:c6}
\end{align}
In particular, if $N$ is odd then $d$ is even.
\end{proposition}

\begin{proof}
Since $G$ is induced, each $v\in U$ has exactly $d_G(v)$ neighbours in $U$ and $d-d_G(v)$ in $W$; summing gives the number of $U$--$W$ edges. Summing degrees over $W$ gives $td=2|E(X[W])|+\sigma_d$. Now \eqref{eq:c2} holds because $X$ is simple and contains $G$; \eqref{eq:c3} because $|E(X[W])|\ge0$; \eqref{eq:c4} because a vertex of degree $\delta$ in $G$ has $d-\delta$ distinct neighbours in $W$; \eqref{eq:c5} because $2|E(X[W])|\le t(t-1)$; and \eqref{eq:c6} by the handshake lemma.
\end{proof}

Conditions \eqref{eq:c3}--\eqref{eq:c6} are those of \citet{EK63}, who showed that for $d=\Delta$ they determine the least order of a $\Delta$-regular graph containing $G$ as an induced subgraph. They are listed, together with the condition for regular graphs mentioned in the introduction, in \citet[Lemma~2.5]{FB26c}. For abelian Cayley hosts they are necessary but not sufficient.

\begin{definition}\label{def:EK}
For an integer $d_0\ge\Delta$, let $\EK(G;d_0)$ be the least $N\ge n$ for which some integer $d$ with $d_0\le d\le N-1$ satisfies \eqref{eq:c3}--\eqref{eq:c6} with $t=N-n$. Thus $\EK(G;\Delta)$ is the floor $\EK(G)$ of \citet[Definition~2.6]{FB26c}, and $\eta(G)\ge\EK(G;d_0)$ whenever every induced embedding of $G$ has degree at least $d_0$.
\end{definition}

\begin{corollary}\label{cor:EKbound}
For every integer $d_0\ge\Delta$,
\[
\EK(G;d_0)\ \ge\ \Bigl\lceil\max\Bigl\{2n-\frac{2m}{d_0},\ n+d_0-\delta\Bigr\}\Bigr\rceil.
\]
\end{corollary}

\begin{proof}
Condition \eqref{eq:c3} is equivalent to $N\ge2n-2m/d$ and \eqref{eq:c4} to $N\ge n+d-\delta$. Both right-hand sides increase with $d$, so they are at least their values at $d_0$.
\end{proof}

\section{The weighted decomposition number}\label{sec:wd}

\begin{definition}\label{def:wd}
Let $G$ be a graph, not necessarily connected. A \emph{weighted $(2,1)$-decomposition} of $G$ partitions $E(G)$ into sets $A_1,\dots,A_a$ and $B_1,\dots,B_b$, some possibly empty, with $\Delta(G[A_i])\le2$ and $\Delta(G[B_j])\le1$. Its \emph{weight} is $2a+b$, and $\wdd(G)$ is the least weight; $\wdd(G)=0$ if $E(G)=\emptyset$.
\end{definition}

\begin{proposition}\label{prop:wdbounds}
For every graph $G$, $\Delta\le\wdd(G)\le\chi'(G)\le\Delta+1$.
\end{proposition}

\begin{proof}
For a decomposition as in Definition~\ref{def:wd} and any vertex $v$, $d_G(v)=\sum_i d_{G[A_i]}(v)+\sum_j d_{G[B_j]}(v)\le2a+b$; take $v$ of degree $\Delta$. A proper edge colouring with $\chi'(G)$ colours is a decomposition into matchings. The last inequality is Vizing's theorem \citep{Vi64}.
\end{proof}

\begin{lemma}[\citealp{Pe1891}]\label{lem:petersen}
Let $F$ be a graph, not necessarily connected, and $k\ge0$ an integer with $\Delta(F)\le2k$. Then $E(F)$ can be partitioned in polynomial time into $k$ sets, each inducing maximum degree at most two.
\end{lemma}

\begin{proof}
For $k=0$ there is nothing to prove. Let $k\ge1$ and let $T$ be the set of odd-degree vertices; $|T|$ is even. Add a new vertex $z$ joined to every vertex of $T$ (if $T\ne\emptyset$). Every vertex of the new graph $F^+$ has even degree, so every component has an Euler tour; orient the edges along these tours. Then every original vertex $v$ has in- and out-degree $\frac12d_{F^+}(v)\le k$, since $d_F(v)\le2k-1$ if $d_F(v)$ is odd. Delete $z$, keeping the orientations. Form the bipartite graph with vertex classes $\{v^+\}$ and $\{v^-\}$ and an edge $u^+v^-$ for each arc $u\to v$. Its maximum degree is at most $k$, so by K\H{o}nig's edge-colouring theorem \citep{Ko16} its edges split into $k$ matchings, found in polynomial time by repeated bipartite matching. A matching uses at most one arc out of and one arc into each vertex, so the corresponding edge set of $F$ has maximum degree at most two. For a modern treatment of this theorem and of the matching results used below, see \citet{LP86}.
\end{proof}

\begin{theorem}\label{thm:wd}
Let $G$ be a graph with $\Delta=\Delta(G)\ge1$, and let $X_\Delta=\{v:d_G(v)=\Delta\}$. Then
\[
\wdd(G)=\begin{cases}\Delta+1, & \text{if $\Delta$ is odd and no matching of $G$ covers $X_\Delta$,}\\ \Delta, & \text{otherwise.}\end{cases}
\]
An optimal weighted $(2,1)$-decomposition can be found in polynomial time.
\end{theorem}

\begin{proof}
By Proposition~\ref{prop:wdbounds}, $\wdd(G)\ge\Delta$. If $\Delta=2k$, Lemma~\ref{lem:petersen} gives $k$ parts of maximum degree at most two, of weight $\Delta$.

Let $\Delta=2k+1$. If a matching $M$ covers $X_\Delta$, then $\Delta(G-M)\le2k$, and Lemma~\ref{lem:petersen} applied to $G-M$, together with the matching $M$, gives weight $2k+1$. Conversely, suppose a decomposition with $a$ parts of maximum degree two and $b$ matchings has weight $2a+b=\Delta$. Then $b$ is odd, so $b\ge1$. For $x\in X_\Delta$,
\[
\Delta=d_G(x)=\sum_{i=1}^a d_{G[A_i]}(x)+\sum_{j=1}^b d_{G[B_j]}(x)\le2a+b=\Delta,
\]
so every summand attains its capacity; in particular each $B_j$ covers $X_\Delta$. Hence $\wdd(G)=\Delta$ exactly when a matching covers $X_\Delta$, and otherwise Lemma~\ref{lem:petersen} with $k+1$ parts gives $\wdd(G)=\Delta+1$.

A matching covering $X_\Delta$ exists exactly when a maximum-weight matching for the edge weights $w(uv)=[u\in X_\Delta]+[v\in X_\Delta]$ has weight $|X_\Delta|$; such a matching is found in polynomial time \citep{Ed65}.
\end{proof}

\begin{corollary}\label{cor:cubic}
A cubic graph has $\wdd=3$ if it has a perfect matching and $\wdd=4$ otherwise. In particular the Petersen graph has chromatic index four and $\wdd=3$.
\end{corollary}

Let $\eall(G)$ be the least order of a finite group, abelian or not, having an undirected Cayley graph that contains $G$ as an induced subgraph; clearly $\eall(G)\le\eta(G)$.

\begin{theorem}\label{thm:nonab}
Let $\Gamma$ be a finite group, written multiplicatively, let $S=S^{-1}\subseteq\Gamma\setminus\{1\}$, and suppose that $G$ is an induced subgraph of $\Cay(\Gamma,S)$, in which $x$ is adjacent to $xs$ for $s\in S$. Then $|S|\ge\wdd(G)$. Consequently $\eta(G)\ge\eall(G)\ge\EK(G;\wdd(G))$, and in particular
\begin{equation}\label{eq:wdbound}
\eta(G)\ \ge\ \eall(G)\ \ge\ \Bigl\lceil\max\Bigl\{2n-\frac{2m}{\wdd(G)},\ n+\wdd(G)-\delta(G)\Bigr\}\Bigr\rceil.
\end{equation}
\end{theorem}

\begin{proof}
An edge $\{x,y\}$ of the host has the unique label $\{s,s^{-1}\}$ with $s=x^{-1}y$, so the classes $\{s,s^{-1}\}$ partition the edges of the host and of the copy of $G$. For $s\in S$ with $s\ne s^{-1}$, the edges $\{x,xs\}$ form one cycle $x,xs,xs^2,\dots,xs^{\operatorname{ord}(s)-1}$ for each orbit of right multiplication by $s$. For $s=s^{-1}$, right multiplication by $s$ is a fixed-point-free involution, and its edges form a perfect matching. As in Lemma~\ref{lem:classes}, the classes $\{s,s^{-1}\}$ restrict to a weighted $(2,1)$-decomposition of $G$ of weight at most $|S|$, so $|S|\ge\wdd(G)$. Every undirected Cayley graph is $|S|$-regular, so Proposition~\ref{prop:EK} applies with $d=|S|\ge\wdd(G)$, and Corollary~\ref{cor:EKbound} gives \eqref{eq:wdbound}.
\end{proof}

\begin{remark}\label{rem:ceiling}
The right-hand side of \eqref{eq:wdbound} is at most $2n-1$: since $G$ is connected, $2m/\wdd(G)\ge2(n-1)/n\ge1$, and $n+\wdd(G)-\delta\le n+\Delta+1-1\le2n-1$. This ceiling does not limit $\eta$: by \citet[Theorem~5.1]{FB26c}, $\eta(K_{1,1,q})=3q$, which exceeds $2n-1=2q+3$ for $q\ge4$.
\end{remark}

\begin{remark}\label{rem:regular}
Let $G$ be regular with $\wdd(G)=\Delta+1$. Then $\Delta$ is odd and $G$ has no perfect matching. A Cayley graph of odd degree has a perfect matching, because a connection set of odd size contains an element $s=s^{-1}$; so $G$ is not a Cayley graph, and \citet[Lemma~2.5(v)]{FB26c}, whose proof uses only that the components of a Cayley graph are Cayley graphs, already gives $d\ge\Delta+1$. Thus for regular graphs Theorem~\ref{thm:nonab} adds nothing to the floor of \citet{FB26c}. For irregular graphs it can, but even then the gain may vanish, because \eqref{eq:c5} and \eqref{eq:c6} may already exclude $d=\Delta$: for the graph $B$ of Section~\ref{sec:abelian}, $\wdd(B)=\Delta+1$ but $\EK(B;\Delta)=\EK(B;\Delta+1)=14$. Example~\ref{ex:Bprime} gives a graph for which the floor rises. For the graphs of Section~\ref{sec:abelian} the rigidity results give stronger bounds on abelian hosts, but they say nothing about non-abelian hosts, to which Theorem~\ref{thm:nonab} applies.
\end{remark}

In the join bound of \citet[Theorem~3.1]{FB26c}, the clique has $c$ vertices; we write $\omega$ for this number, since $c$ denotes a generator class here. For $\omega=1$, it states that $\eta(G)\ge2\max_v\alpha(G[N_G(v)])$; this is the abelian form of \citet[Proposition~7.4]{BS85}. The first term of \eqref{eq:wdbound} dominates it.

\begin{corollary}\label{cor:join}
For every connected graph $G$,
\[
2n-\frac{2m}{\wdd(G)}\ \ge\ 2n-\frac{2m}{\Delta}\ \ge\ \max\{n,2\alpha(G)\}\ \ge\ \max\Bigl\{n,\,2\max_v\alpha(G[N_G(v)])\Bigr\}.
\]
\end{corollary}

\begin{proof}
The first inequality holds because $\wdd(G)\ge\Delta$, and $2m\le n\Delta$ gives $2n-2m/\Delta\ge n$. If $A$ is independent, every edge has an end in $V\setminus A$, so $m\le\Delta(n-|A|)$, that is, $2n-2m/\Delta\ge2|A|$. Finally $\alpha(G[N_G(v)])\le\alpha(G)$.
\end{proof}

For $\omega\ge2$ the join bound is not dominated: for $K_{1,1,q}$ with $q\ge2$ it gives the exact value $3q$, while \eqref{eq:wdbound} gives $2q+1$. Conversely, for the graph $B$ of Section~\ref{sec:abelian}, \eqref{eq:wdbound} gives $14$ and every join bound is at most $6$.

\begin{example}\label{ex:terms}
Neither term of \eqref{eq:wdbound} dominates the other.
\begin{enumerate}
\item[(a)] The double star $D_{q,q}$, two adjacent centres each with $q\ge2$ leaves, has $n=2q+2$, $m=2q+1$, $\Delta=q+1$, $\delta=1$, and $\wdd=\Delta$ because the central edge covers both vertices of maximum degree. The two terms are $\lceil4q+2/(q+1)\rceil=4q+1$ and $3q+2$. Here $\eta(D_{q,q})=5q$ \citep[Corollary~6.5]{FB26c}.
\item[(b)] The graph $K_q^+$ obtained from $K_q$, $q\ge4$, by attaching a pendant vertex has $n=q+1$, $m=\binom q2+1$, $\Delta=q$, $\delta=1$ and $\wdd=q$. The two terms are $\lceil q+3-2/q\rceil=q+3$ and $2q$. For $q=4$ the second is exact: $\eta(K_4^+)=8$ \citep[Proposition~7.1]{FB26c}.
\end{enumerate}
\end{example}

\section{Abelian obstructions}\label{sec:abelian}

The weighted invariant uses only the degree conditions of Lemma~\ref{lem:classes}. In an abelian group, two generators $s,g$ also satisfy $(x+s)+g=(x+g)+s$, and this relation gives much stronger restrictions.

\subsection{Rigidity of the subdivided \texorpdfstring{$K_4$}{K4}}

\begin{definition}\label{def:J}
Let $J$ be the graph obtained from $K_4$ on $\{\jv a,\jv b,\jv c,\jv d\}$ by subdividing the edge $\jv c\jv d$ with a new vertex $\jv s$:
\[
V(J)=\{\jv a,\jv b,\jv c,\jv d,\jv s\},\qquad E(J)=\{\jv a\jv b,\jv a\jv c,\jv a\jv d,\jv b\jv c,\jv b\jv d,\jv c\jv s,\jv d\jv s\}.
\]
The vertices $\jv c$ and $\jv d$ are non-adjacent and have the three common neighbours $\jv a,\jv b,\jv s$, of which exactly $\jv a$ and $\jv b$ are adjacent. The vertices $\jv a,\jv b,\jv c,\jv d$ induce a diamond, that is, $K_4$ minus the edge $\jv c\jv d$. The vertices of $J$, and of the graphs built from copies of $J$, are written in sans-serif type to distinguish them from group elements, generator classes and degrees.
\end{definition}

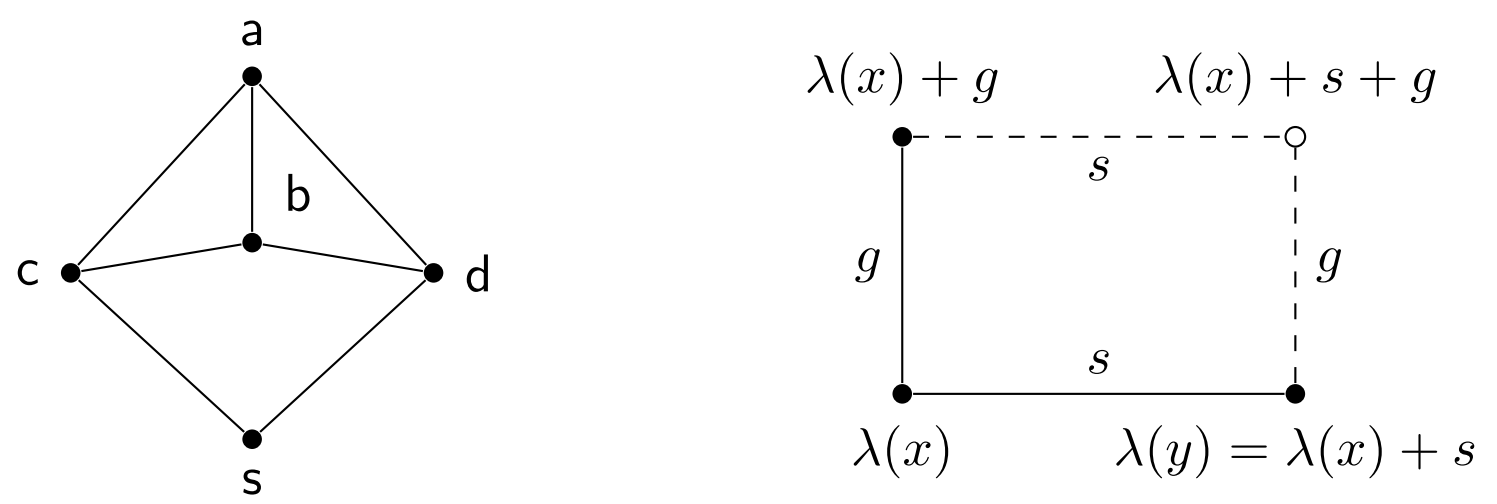
\begin{figure}[ht]
\centering
\begin{tikzpicture}[v/.style={circle,fill,inner sep=1.3pt},font=\small]
\begin{scope}
\node[v,label=above:$\jv a$] (a) at (0,1.3) {};
\node[v,label=left:$\jv c$] (c) at (-1.2,0) {};
\node[v,label={[label distance=1pt]above right:$\jv b$}] (b) at (0,0.2) {};
\node[v,label=right:$\jv d$] (d) at (1.2,0) {};
\node[v,label=below:$\jv s$] (s) at (0,-1.1) {};
\draw (a)--(b) (a)--(c) (a)--(d) (b)--(c) (b)--(d) (c)--(s) (d)--(s);
\node at (0,-1.9) {(a)};
\end{scope}
\begin{scope}[xshift=5.6cm]
\node[v,label=below:$\lambda(x)$] (u) at (-1.3,-0.8) {};
\node[v,label=below:{$\lambda(y)=\lambda(x)+s$}] (y) at (1.3,-0.8) {};
\node[v,label=above:{$\lambda(x)+g$}] (t) at (-1.3,0.9) {};
\node[draw,circle,inner sep=1.3pt,label=above:{$\lambda(x)+s+g$}] (st) at (1.3,0.9) {};
\draw (u)--node[above]{$s$}(y) (u)--node[left]{$g$}(t);
\draw[dashed] (t)--node[below]{$s$}(st) (y)--node[right]{$g$}(st);
\node at (0,-1.9) {(b)};
\end{scope}
\end{tikzpicture}
\caption{(a) The graph $J$: $K_4$ on $\{\jv a,\jv b,\jv c,\jv d\}$ with the edge $\jv c\jv d$ subdivided by $\jv s$. (b) The commuting square of Proposition~\ref{prop:square}: if $G$ has no four-cycle, the corner $\lambda(x)+s+g$ lies outside the image and is a host neighbour of $\lambda(y)$.}
\label{fig:J}
\end{figure}

\begin{lemma}\label{lem:C7}
For $r\ge5$ let $C_r^2=\Cay(\Z_r,\{\pm1,\pm2\})$. Then $C_r^2$ contains an induced copy of $J$ if and only if $r=7$.
\end{lemma}

\begin{proof}
In $J$ the non-adjacent vertices $\jv c,\jv d$ have three common neighbours inducing exactly one edge. For $r=5$, $C_5^2=K_5$ has no non-adjacent pair. For $r=6$, $C_6^2$ is $K_6$ minus the perfect matching $\{x,x+3\}$; the common neighbourhood of a non-adjacent pair consists of the other four vertices, which induce a $4$-cycle, and every three vertices of a $4$-cycle induce two edges. For $r\ge8$, translate a non-adjacent pair to $\{0,h\}$ with $3\le h\le\lfloor r/2\rfloor$. Comparing $N(0)=\{\pm1,\pm2\}$ with $N(h)=\{h\pm1,h\pm2\}$ gives at most two common neighbours: $1,2$ for $h=3$; $2$, and also $-2$ when $r=8$, for $h=4$; none for $h\ge5$. For $r=7$, the labels $\jv c=-1$, $\jv d=2$, $\jv a=0$, $\jv b=1$, $\jv s=4$ work: $\jv c-\jv d=\pm3$ is not a generator; each of $\jv a,\jv b,\jv s$ differs from $\jv c$ and from $\jv d$ by $\pm1$ or $\pm2$; $\jv a-\jv b=-1$ is a generator, while $\jv a-\jv s$ and $\jv b-\jv s$ are $\pm3$.
\end{proof}

\begin{theorem}\label{thm:deg4}
Let $X=\Cay(\Gamma,S)$ be an abelian Cayley graph with $|S|=4$. If $X$ contains an induced copy of $J$, then $\langle S\rangle$ is cyclic of order seven and $\Cay(\langle S\rangle,S)\cong C_7^2$; in particular every component of $X$ has seven vertices. Conversely, $C_7^2$ contains an induced copy of $J$.
\end{theorem}

\begin{proof}
An induced copy of $J$ contains an induced diamond: the triangles $\jv a\jv b\jv c$ and $\jv a\jv b\jv d$ share the edge $\jv a\jv b$, and their third vertices $\jv c,\jv d$ are non-adjacent. We show that this configuration forces $\Cay(\langle S\rangle,S)\cong C_7^2$. Since $S=-S$ and $|S|=4$, there are three cases.

\emph{Four involutions.} Translate an edge to $\{0,g\}$ with $g\in S$. A third vertex of a triangle on this edge is an element $h\in S\setminus\{g\}$ with $g-h=g+h\in S\setminus\{g\}$. Such elements come in pairs $\{h,g+h\}$, and $S\setminus\{g\}$ has three elements, so the edge lies in no triangle or in exactly two, whose third vertices $h$ and $g+h$ differ by $g\in S$ and are adjacent. So no diamond occurs. (Two triangles on an edge do occur: in $\Cay(\Z_2^3,\{e_1,e_2,e_1+e_2,e_3\})\cong K_4\,\square\,K_2$ the edge $\{0,e_1\}$ lies in two triangles.)

\emph{Two involutions and one inverse pair.} Write $S=\{j_1,j_2,g,-g\}$ with $j_1,j_2$ involutions and $2g\ne0$. The three steps of a triangle lie in $S$ and sum to zero. They are not all involutions, since they would be distinct and $S$ has only two; and exactly two cannot be involutions, since the sum of two involutions is zero or an involution. So a triangle either uses only $\pm g$, which requires $\operatorname{ord}(g)=3$, or uses one involution and two steps $\pm g$, which requires $2g\in\{j_1,j_2\}$. If $2g\notin\{j_1,j_2\}$, every triangle uses only $\pm g$, each edge lies in at most one triangle, and no diamond occurs; so we may assume $2g=j_1$. Then $\operatorname{ord}(g)=4$, $j_2\notin\langle g\rangle$, $\langle S\rangle\cong\Z_4\times\Z_2$ with $S=\{(\pm1,0),(2,0),(0,1)\}$, and $\Cay(\langle S\rangle,S)\cong K_4\,\square\,K_2$. Every triangle of $K_4\,\square\,K_2$ lies in one copy of $K_4$, so two triangles sharing an edge have adjacent third vertices, and no diamond occurs.

\emph{Two inverse pairs.} Write $S=\{\pm g,\pm h\}$ with $2g,2h\ne0$. We claim that an edge in two triangles forces $h=\pm2g$ or $g=\pm2h$. Translate the edge to $\{0,g\}$, after renaming. A third vertex lies in $\{-g,h,-h\}$ and is adjacent to $g$. If it is $-g$, then $-2g\in S$, which without a relation $h=\pm2g$ means $-2g=g$, that is $\operatorname{ord}(g)=3$, and gives one triangle. If it is $\pm h$, then $\pm h-g\in S$, which forces $h=\pm2g$ or $g=\pm2h$. Indeed, the eight cases are $h-g\in\{g,-g,h,-h\}$ and $-h-g\in\{g,-g,h,-h\}$; they give $h=2g$, $h=0$, $g=0$, $g=2h$ and $h=-2g$, $h=0$, $g=-2h$, $g=0$ respectively, and $g,h\ne0$. So without the relation every edge lies in at most one triangle. Hence, after renaming, $S=\{\pm g,\pm2g\}$, $\langle S\rangle=\langle g\rangle$, and $r=\operatorname{ord}(g)\ge5$ because the four elements are distinct and not involutions. Then $\Cay(\langle S\rangle,S)\cong C_r^2$. The connected copy of $J$ lies in one component (Lemma~\ref{lem:connected}), so $C_r^2$ contains an induced $J$, and Lemma~\ref{lem:C7} gives $r=7$. All components are cosets of $\langle S\rangle$ by Lemma~\ref{lem:connected}. The converse is Lemma~\ref{lem:C7}.
\end{proof}

\begin{corollary}\label{cor:deg4}
If $G$ is connected, contains an induced copy of $J$ and has at least eight vertices, then $G$ has no induced embedding of degree four.
\end{corollary}

\begin{proof}
By Lemma~\ref{lem:connected} the image of $G$ would lie in one component of the host, which has seven vertices by Theorem~\ref{thm:deg4}.
\end{proof}

For degree five we use a finite symbolic computation instead of a case analysis by hand.

\begin{theorem}\label{thm:deg5}
Let $X=\Cay(\Gamma,S)$ be an abelian Cayley graph with $|S|\le5$ that contains an induced copy of $J$. Then either $|S|=4$ and $|\langle S\rangle|=7$, or $|S|=5$ and $|\langle S\rangle|\in\{8,10,14\}$. All three orders occur for $|S|=5$. For $|S|=6$ the order of $\langle S\rangle$ is unbounded.
\end{theorem}

\begin{proof}
Since $J$ has vertices of degree three, $|S|\ge3$. Write $S=\{\pm x_1,\dots,\pm x_k,i_1,\dots,i_j\}$ with $2k+j=|S|$, where $i_1,\dots,i_j$ are the involutions in $S$. Translate an induced copy of $J$ so that $\jv a$ has label $0$, and let $\sigma_1,\dots,\sigma_7\in S$ be the differences along its seven edges:
\[
\jv b=\sigma_1,\quad \jv c=\sigma_2,\quad \jv d=\sigma_3,\quad \jv c-\jv b=\sigma_4,\quad \jv d-\jv b=\sigma_5,\quad \jv s-\jv c=\sigma_6,\quad \jv s-\jv d=\sigma_7,
\]
where vertices stand for their labels. Let $\mathcal P$ be the abelian group generated by symbols $x_1,\dots,x_k,i_1,\dots,i_j$ subject to $2i_t=0$ and to the three relations $\sigma_2-\sigma_1=\sigma_4$, $\sigma_3-\sigma_1=\sigma_5$ and $\sigma_2+\sigma_6=\sigma_3+\sigma_7$, each $\sigma_r$ being read as the corresponding symbol. Sending each symbol to the element it names defines a homomorphism from $\mathcal P$ onto $\langle S\rangle$. The other requirements are inequalities: the $|S|$ symbols are distinct and non-zero, and $2x_t\ne0$; the five vertices are distinct; and $\jv c-\jv d$, $\jv s$ and $\jv s-\jv b$ lie outside $S\cup\{0\}$ (the non-edges $\jv c\jv d$, $\jv a\jv s$ and $\jv b\jv s$). An inequality that fails in $\mathcal P$ fails in every quotient of $\mathcal P$.

For each $|S|\in\{3,4,5\}$, each split $2k+j=|S|$ and each of the $|S|^7$ choices of $\sigma_1,\dots,\sigma_7$, we computed the Hermite normal form of the relation lattice. It gives the order of $\mathcal P$, or shows that $\mathcal P$ is infinite, and decides every inequality. For $|S|=3$ no choice satisfies all inequalities in $\mathcal P$. For $|S|=4$ every choice that does has $|\mathcal P|=7$, and for $|S|=5$ every such choice has $|\mathcal P|\in\{8,10,14\}$; no such $\mathcal P$ is infinite. Hence $|\langle S\rangle|$ divides $7$ in the first case and one of $8,10,14$ in the second. In the second case the divisors that are at most $7$ are $1,2,4,5$ and $7$, and none of them occurs: a group of order $7$ has no involution, so its connection sets have even size, and in a group of order at most five the copy of $J$ would be the whole Cayley graph, which is regular, whereas $J$ is not. This proves the first statement. The connection sets $\{1,2,4,6,7\}\subseteq\Z_8$, $\{1,2,5,8,9\}\subseteq\Z_{10}$ and $\{2,4,7,10,12\}\subseteq\Z_{14}$ give induced copies of $J$ with $\langle S\rangle$ of order $8$, $10$ and $14$. For $|S|=6$, the graph $\Cay(\Z_7\times\Z_m,\{(\pm1,0),(\pm2,0),(0,\pm1)\})\cong C_7^2\,\square\,C_m$, $m\ge3$, contains the induced copy of $J$ of Lemma~\ref{lem:C7} in the layer $\Z_7\times\{0\}$.
\end{proof}

\begin{remark}[Reproducible enumeration protocol]\label{rem:protocol}
The finite assertion in Theorem~\ref{thm:deg5} can be reproduced with the script \texttt{theorem45\_lattice.py}, and checked independently by the direct connection-set test \texttt{brute45\_groups.py}; both are in the data deposit described at the end of the paper. For each $d=3,4,5$ the algorithm enumerates every split $2k+j=d$, constructs the formal signed generator list, and loops over all $d^7$ assignments of the seven edge differences. It forms the three edge-consistency relations and the $j$ involution relations, reduces their integer row lattice to Hermite normal form, and tests every required nonzero, distinctness and nonedge inequality modulo the row lattice. For a consistent presentation it records the rank and product of the nonzero pivots. The exact reproducible counts are
\begin{center}
\begin{tabular}{rrrl}
$d$ & assignments & valid presentations & orders of valid presentations\\\midrule
3 & 4374 & 0 & $\emptyset$\\
4 & 49152 & 16 & 7 (16)\\
5 & 234375 & 448 & 8 (400), 10 (32), 14 (16).
\end{tabular}
\end{center}
Thus the complete search contains $287\,901$ assignments. Its finite enumeration, integer arithmetic and inequality tests can be checked without a SAT solver. For $|S|=4$ the computation reproduces Theorem~\ref{thm:deg4}, which serves as a check of the method.
\end{remark}

\begin{corollary}\label{cor:deg6}
Let $G$ be connected and contain an induced copy of $J$. If $\eta(G)\ge15$, in particular if $G$ has at least fifteen vertices, then every induced embedding of $G$ has degree at least six.
\end{corollary}

\begin{proof}
An induced embedding of degree at most five has degree four or five by Theorem~\ref{thm:deg5}. By Lemma~\ref{lem:connected} it can be translated into $\Cay(\langle S\rangle,S)$, which has at most $14$ vertices by Theorem~\ref{thm:deg5}. This contradicts $\eta(G)\ge15$.
\end{proof}

\subsection{Commuting squares}
If $s,g\in S$ and $g\ne\pm s$, then $0,s,s+g,g$ is a four-cycle of $\Cay(\Gamma,S)$; this is why abelian Cayley graphs of degree at least three have girth at most four. For graphs without four-cycles, the same square forces host neighbours outside the image.

\begin{proposition}\label{prop:square}
Let $G=(V,E)$ be a graph without four-cycles (induced or not), and let $(\Gamma,S,\lambda)$ be an induced embedding of degree $d$. For an ordered pair $(x,y)$ with $xy\in E$ put $s_{xy}=\lambda(y)-\lambda(x)$, $A_x=\{\lambda(z)-\lambda(x):z\in N_G(x)\}\subseteq S$, and let $\epsilon_{xy}=1$ if $-s_{xy}\in A_x$ and $-s_{xy}\ne s_{xy}$, and $\epsilon_{xy}=0$ otherwise. Then
\begin{equation}\label{eq:square}
d-d_G(y)\ge d_G(x)-1-\epsilon_{xy},\qquad\text{hence}\qquad d\ge d_G(x)+d_G(y)-1-\epsilon_{xy}.
\end{equation}
\end{proposition}

\begin{proof}
Put $u=\lambda(x)$ and $s=s_{xy}$, so $\lambda(y)=u+s$. Since $\lambda$ is injective, $|A_x|=d_G(x)$ and $s\in A_x$. Let $T=A_x\setminus\{s,-s\}$; then $|T|=d_G(x)-1-\epsilon_{xy}$. For $g\in T$ the elements $u$, $u+s$, $u+s+g$, $u+g$ are distinct and form a four-cycle of the host (Figure~\ref{fig:J}(b)). Three of them are labels: $u+s$ of $y$, and $u+g$ of the neighbour of $x$ reached by $g$. If $u+s+g$ were a label too, inducedness would make the four-cycle a four-cycle of $G$. So $u+s+g\notin\lambda(V)$, and it is a host neighbour of $u+s=\lambda(y)$. Distinct $g$ give distinct $u+s+g$, so $y$ has at least $|T|$ host neighbours outside $\lambda(V)$, and it has exactly $d-d_G(y)$ of them.
\end{proof}

\begin{corollary}\label{cor:C4free}
Let $G$ be a connected $k$-regular graph without four-cycles, $k\ge1$, and let $b_k=2k-1$ for odd $k$ and $b_k=2k-2$ for even $k$. Every induced embedding of $G$ has degree $d\ge\max\{k,b_k\}$, and
\[
\eta(G)\ \ge\ \EK(G;\max\{k,b_k\})\ \ge\ \Bigl\lceil 2n-\frac{kn}{\max\{k,b_k\}}\Bigr\rceil.
\]
For $k=3$ the degree bound is attained (Proposition~\ref{prop:desargues}).
\end{corollary}

\begin{proof}
Proposition~\ref{prop:square} gives $d\ge2k-1-\epsilon_{xy}\ge2k-2$ for every edge. If $k$ is odd, the map $s\mapsto-s$ cannot pair off the $k$ elements of $A_x$, so some $s\in A_x$ has $-s\notin A_x$ or $-s=s$; for the corresponding edge $\epsilon_{xy}=0$ and $d\ge2k-1$. Corollary~\ref{cor:EKbound} with $2m=kn$ gives the bound on $\eta$.
\end{proof}

\begin{corollary}\label{cor:girth5}
Let $G$ be a connected cubic graph of girth at least five on $n$ vertices. Then every induced embedding of $G$ has degree at least five, and
\[
\eta(G)\ \ge\ 2\Bigl\lceil\frac{7n}{10}\Bigr\rceil.
\]
In particular $\eta(P)\ge14$ for the Petersen graph $P$.
\end{corollary}

\begin{proof}
Here $n$ is even, and $n\ge10$ because a vertex, its three neighbours and the six vertices at distance two are distinct. By Corollary~\ref{cor:C4free}, $d\ge5$. If $d=5$, then the host order $N$ is even by Proposition~\ref{prop:EK}, and $N\ge2n-3n/5=7n/5$, so $N\ge2\lceil7n/10\rceil$. If $d\ge6$, then $N\ge2n-3n/6=3n/2$, and $2\lceil7n/10\rceil\le3n/2$ for every even $n\ge10$: for $n=2r$ check $r=5,\dots,9$ directly, and for $r\ge10$ use $2\lceil7r/5\rceil<14r/5+2\le3r$.
\end{proof}

For the Petersen graph, Corollary~\ref{cor:girth5} improves the floor $13$ that the Erd\H{o}s--Kelly conditions give with $d\ge\Delta+1$ \citep[Section~2]{FB26c}.

\begin{proposition}\label{prop:desargues}
Let $Q_5=\Cay(\F_2^5,\{e_1,\dots,e_5\})$, and let $U$ be the set of vectors of Hamming weight two or three. Then $Q_5[U]$ is the Desargues graph $D$, which is cubic of girth six, so $D$ has an induced embedding of order $32$ and degree five. Hence the degree bound of Corollary~\ref{cor:C4free} is attained for $k=3$, and $28\le\eta(D)\le32$.
\end{proposition}

\begin{proof}
Identify vectors with subsets of $\{1,\dots,5\}$. Two sets in $U$ are adjacent in $Q_5$ exactly when one is a $2$-set contained in the other, a $3$-set. Every $2$-set lies in three $3$-sets and every $3$-set contains three $2$-sets, so $Q_5[U]$ is cubic and bipartite. Two distinct $2$-sets lie in at most one common $3$-set, and two distinct $3$-sets contain at most one common $2$-set, so there is no four-cycle; $\{1,2\},\{1,2,3\},\{2,3\},\{2,3,4\},\{2,4\},\{1,2,4\}$ is a six-cycle. Replacing each $3$-set by its complement identifies $Q_5[U]$ with the bipartite Kneser graph $H(5,2)$, which is the Desargues graph. Corollary~\ref{cor:girth5} with $n=20$ gives $\eta(D)\ge28$.
\end{proof}

\subsection{Two extremal classifications}

\begin{lemma}\label{lem:33332}
Every graph on five vertices with degree sequence $(3,3,3,3,2)$ is isomorphic to $J$, with the vertex of degree two corresponding to $\jv s$.
\end{lemma}

\begin{proof}
The complement has degree sequence $(1,1,1,1,2)$ and three edges. Its vertex of degree two is adjacent to two vertices of degree one, and the remaining two vertices of degree one are adjacent to each other, so the complement is $P_3\cup K_2$. The complement of $P_3\cup K_2$, with the middle vertex of $P_3$ as $\jv s$ and the edge $K_2$ as $\jv c\jv d$, is $J$.
\end{proof}

\begin{lemma}[Relative Tutte criterion]\label{lem:tutte}
Let $G$ be a graph and $X\subseteq V(G)$. For $U\subseteq V(G)$ let $o_X(G-U)$ be the number of components of $G-U$ of odd order contained in $X$. Then $G$ has a matching covering $X$ if and only if $o_X(G-U)\le|U|$ for every $U\subseteq V(G)$.
\end{lemma}

\begin{proof}
This subset-saturating form of Tutte's theorem is classical \citep[see, e.g.,][]{LP86}; we include the short proof for self-containment. If a matching covers $X$, every odd component of $G-U$ contained in $X$ has a vertex matched into $U$, and distinct components use distinct vertices of $U$. Conversely, suppose that no matching covers $X$. Add a set $Z$ of new vertices with $|Z|\ge|V(G)|$ and $|Z|\equiv|V(G)|\pmod2$, all edges inside $Z$, and all edges between $Z$ and $V(G)\setminus X$. A perfect matching of the new graph $G^*$ would cover $X$ with edges of $G$, so $G^*$ has none. By Tutte's theorem \citep{Tu47} there is $U^*$ with $o(G^*-U^*)>|U^*|$, and by parity $o(G^*-U^*)\ge|U^*|+2$. If $Z\subseteq U^*$, then $o(G^*-U^*)\le|V(G)|\le|U^*|$, which is impossible. So $Z\not\subseteq U^*$, and the vertices of $Z\setminus U^*$ and of $V(G)\setminus(X\cup U^*)$ lie in one component of $G^*-U^*$. Every other component is contained in $X$ and is a component of $G-U$, where $U=U^*\cap V(G)$. Hence $o_X(G-U)\ge o(G^*-U^*)-1\ge|U^*|+1>|U|$.
\end{proof}

\begin{figure}[ht]
\centering
\begin{tikzpicture}[v/.style={circle,fill,inner sep=1.3pt},font=\small,scale=0.95]
\node[v,label=left:$\jv a_2$] (a2) at (-4.6,0) {};
\node[v,label=above:$\jv c_2$] (c2) at (-3.6,0.9) {};
\node[v,label=below:$\jv d_2$] (d2) at (-3.6,-0.9) {};
\node[v,label={[label distance=0pt]right:$\jv b_2$}] (b2) at (-3.7,0) {};
\node[v,label=above:$\jv s_2$] (s2) at (-2.4,0) {};
\draw (a2)--(b2) (a2)--(c2) (a2)--(d2) (b2)--(c2) (b2)--(d2) (c2)--(s2) (d2)--(s2);
\node[v,label=above:$\jv z$] (z) at (-1.4,0) {};
\node[v,label=above:$\jv s_1$] (s1) at (-0.4,0) {};
\node[v,label=above:$\jv c_1$] (c1) at (0.8,0.9) {};
\node[v,label=below:$\jv d_1$] (d1) at (0.8,-0.9) {};
\node[v,label={[label distance=0pt]left:$\jv b_1$}] (b1) at (0.9,0) {};
\node[v,label=right:$\jv a_1$] (a1) at (1.8,0) {};
\draw (a1)--(b1) (a1)--(c1) (a1)--(d1) (b1)--(c1) (b1)--(d1) (c1)--(s1) (d1)--(s1);
\draw (s2)--(z)--(s1);
\end{tikzpicture}
\caption{The graph $B$: two copies of $J$ attached at subdivision vertices $\jv s_1,\jv s_2$ to a degree-two vertex $\jv z$. All vertices are labelled; the labels follow the edge list in~\eqref{eq:EJ}.}
\label{fig:B}
\end{figure}
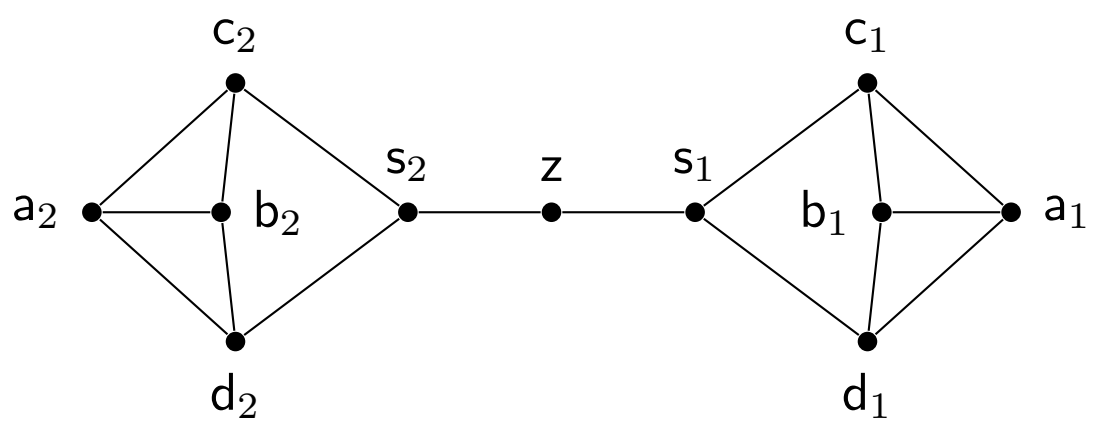

\begin{definition}\label{def:B}
Let $J_1,J_2$ be disjoint copies of $J$, with vertices $\jv a_i,\jv b_i,\jv c_i,\jv d_i,\jv s_i$ and
\begin{equation}\label{eq:EJ}
E(J_i)=\{\jv a_i\jv b_i,\ \jv a_i\jv c_i,\ \jv a_i\jv d_i,\ \jv b_i\jv c_i,\ \jv b_i\jv d_i,\ \jv c_i\jv s_i,\ \jv d_i\jv s_i\}.
\end{equation}
The graph $B$ is obtained by adding a vertex $\jv z$ adjacent to $\jv s_1$ and $\jv s_2$ (Figure~\ref{fig:B}).
\end{definition}

\begin{proposition}\label{prop:B}
The graph $B$ is connected, with $11$ vertices, $16$ edges, $\Delta=3$, $\delta=2$ and $\wdd(B)=4$. Every induced embedding of $B$ has degree at least six, and $\eta(B)\ge17$.
\end{proposition}

\begin{proof}
The ten vertices other than $\jv z$ have degree three. A matching covering them would cover each $J_i$, which has odd order and only the edge $\jv s_i\jv z$ leaving it; so it would contain both $\jv s_1\jv z$ and $\jv s_2\jv z$, which is impossible. Theorem~\ref{thm:wd} gives $\wdd(B)=4$. Since $B$ contains induced copies of $J$, Theorem~\ref{thm:deg5} excludes $d=3$ and Corollary~\ref{cor:deg4} excludes $d=4$. For $d\ge5$, Corollary~\ref{cor:EKbound} gives $N\ge22-32/5$, so $\eta(B)\ge16$. Now Corollary~\ref{cor:deg6} applies and gives $d\ge6$, and then $N\ge22-32/6$, so $\eta(B)\ge17$. Theorem~\ref{thm:nonab} also excludes $d=3$, without computation; unlike the results of this section it applies to hosts in non-abelian groups, where it gives $\eall(B)\ge\EK(B;4)=14$.
\end{proof}

\begin{theorem}\label{thm:2D5}
Let $G$ be connected with $\wdd(G)=\Delta+1$. Then $\Delta\ge3$ is odd and $|V(G)|\ge2\Delta+5\ge11$. Up to isomorphism, $B$ is the only such graph on eleven vertices.
\end{theorem}

\begin{proof}
By Theorem~\ref{thm:wd}, $\Delta$ is odd and no matching covers $X=X_\Delta$. By Lemma~\ref{lem:tutte} there is $U$ with $q=o_X(G-U)\ge|U|+1$; put $u=|U|$ and call the $q$ odd components of $G-U$ contained in $X$ \emph{full}. Every vertex of a full component $C$ has degree $\Delta$, so the number $e(C,U)$ of edges from $C$ to $U$ is $\Delta|C|-2|E(G[C])|$, which is odd. If $u=0$, then $G$ is a full component, a $\Delta$-regular graph of odd order with $\Delta$ odd, which is impossible; so $u\ge1$.

If $|C|\le\Delta$, each vertex of $C$ has at most $|C|-1$ neighbours in $C$, so $e(C,U)\ge|C|(\Delta-|C|+1)\ge\Delta$. Let $\ell$ be the number of full components with $e(C,U)<\Delta$; each has odd order greater than $\Delta$, hence at least $\Delta+2$ vertices. Counting edges at $U$,
\[
\Delta u\ \ge\ \sum_{C\text{ full}}e(C,U)\ \ge\ \ell+\Delta(q-\ell),\qquad\text{so}\qquad \ell(\Delta-1)\ge\Delta(q-u)\ge\Delta,
\]
which is impossible for $\Delta=1$ and gives $\ell\ge2$ for $\Delta\ge3$. Therefore $|V(G)|\ge u+\ell(\Delta+2)+(q-\ell)\ge1+2(\Delta+2)=2\Delta+5$.

If $|V(G)|=11$, then $\Delta=3$, $u=1$, $\ell=q=2$, the two full components have five vertices each, and there are no other vertices. Each full component sends an odd number of edges, fewer than three, to the vertex $w$ of $U$, that is, exactly one. So each full component has degree sequence $(3,3,3,3,2)$ and is a copy of $J$ attached to $w$ at its vertex of degree two (Lemma~\ref{lem:33332}). Hence $G\cong B$.
\end{proof}

\begin{theorem}[Sharpness for every odd maximum degree]\label{thm:sharp}
For every odd integer $\Delta\ge3$ there is a connected simple graph $G_\Delta$ on exactly $2\Delta+5$ vertices such that $\Delta(G_\Delta)=\Delta$ and $\wdd(G_\Delta)=\Delta+1$. Consequently the order lower bound of Theorem~\ref{thm:2D5} is sharp for every allowable value of $\Delta$, not merely for $\Delta=3$.
\end{theorem}

\begin{proof}
Put $r=\Delta+2$, which is odd. Form a graph $Q$ on $r$ vertices as the disjoint union of a path on three vertices and $(\Delta-1)/2$ copies of $K_2$. Let $J_\Delta=\overline Q$ be its complement. In $Q$ exactly one vertex, the middle of the path, has degree two; all other vertices have degree one. In $J_\Delta$ the former therefore has degree $r-1-2=\Delta-1$, and all other vertices have degree $r-1-1=\Delta$. Moreover $J_\Delta$ is connected: its minimum degree is $\Delta-1=r-3>(r-1)/2$ for $r\ge7$, which rules out a disconnected graph. For $r=5$ the graph is the connected subdivided-edge $K_4$. Take two disjoint copies of $J_\Delta$ and add one vertex $\jv z$ adjacent exactly to their two vertices of degree $\Delta-1$. The result $G_\Delta$ is connected and has $2r+1=2\Delta+5$ vertices. All vertices except $\jv z$ have degree $\Delta$ and $\jv z$ has degree two, so the maximum degree is $\Delta$. Let $X_\Delta=V(G_\Delta)\setminus\{\jv z\}$. A matching covering $X_\Delta$ would have to match a vertex of each odd-order copy of $J_\Delta$ outside that copy. The only external edge from each copy runs to $\jv z$, and a matching cannot use both edges incident with $\jv z$. Thus no matching covers $X_\Delta$. The characterization in Theorem~\ref{thm:wd} yields $\wdd(G_\Delta)=\Delta+1$. When $\Delta=3$, each $J_\Delta$ is the subdivided-edge $K_4$, and $G_3=B$.
\end{proof}

\begin{example}\label{ex:Bprime}
Let $B'$ be obtained from $B$ by attaching a path $\jv z\,\ell_1\ell_2$ at $\jv z$. Then $n=13$, $m=18$, $\Delta=3$ and $\delta=1$. Now $\jv z$ has degree three, and $J_1,J_2$ are odd components of $B'-\jv z$ contained in $X_\Delta$, so $\wdd(B')=4$ by Lemma~\ref{lem:tutte} and Theorem~\ref{thm:wd}. The Erd\H{o}s--Kelly floor is $\EK(B';3)=16$, attained with $d=3$ and $t=3$, while $\EK(B';4)=17$. With Corollary~\ref{cor:deg4}, every induced embedding of $B'$ has degree at least five, so $\eta(B')\ge\EK(B';5)=20$; then Corollary~\ref{cor:deg6} gives degree at least six, with the same bound $\EK(B';6)=20$. The true value is $\eta(B')=30$ (Remark~\ref{rem:Bprime}).
\end{example}

\begin{figure}[ht]
\centering
\begin{tikzpicture}[v/.style={circle,fill,inner sep=1.3pt},font=\scriptsize,scale=0.95]
\node[v,label={[label distance=1pt]270:$\jv z$}] (z) at (0,0) {};
\foreach \t/\i in {90/1,210/2,330/3}{
  \pgfmathsetmacro{\tp}{\t+90}\pgfmathsetmacro{\tm}{\t-90}
  \pgfmathsetmacro{\tb}{\t+40}\pgfmathsetmacro{\ts}{\t-35}
  \node[v,label={[label distance=0pt]\ts:$\jv s_{\i}$}] (s\i) at ($(0,0)+(\t:0.9)$) {};
  \node[v,label={[label distance=0pt]\tp:$\jv c_{\i}$}] (c\i) at ($(\t:1.6)+(\tp:0.75)$) {};
  \node[v,label={[label distance=0pt]\tm:$\jv d_{\i}$}] (d\i) at ($(\t:1.6)+(\tm:0.75)$) {};
  \node[v,label={[label distance=-1pt]\tb:$\jv b_{\i}$}] (b\i) at (\t:2.0) {};
  \node[v,label={[label distance=0pt]\t:$\jv a_{\i}$}] (a\i) at (\t:2.85) {};
  \draw (z)--(s\i) (s\i)--(c\i) (s\i)--(d\i) (c\i)--(b\i) (d\i)--(b\i) (c\i)--(a\i) (d\i)--(a\i) (a\i)--(b\i);
}
\end{tikzpicture}
\caption{The cubic graph $H$: three copies of $J$ attached at their subdivision vertices to a common vertex $\jv z$. All vertices are labelled as in~\eqref{eq:EJ}.}
\label{fig:H}
\end{figure}
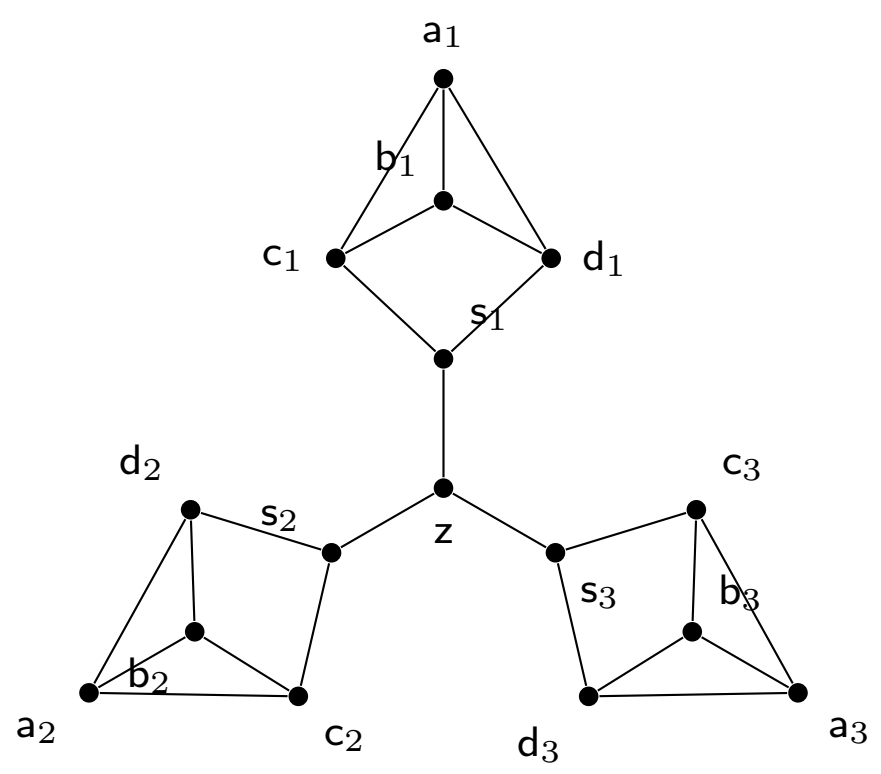

For a tree $T$ with $p\ge1$ vertices and $\Delta(T)\le3$, let $G_T$ be obtained by attaching to each vertex $x$ of $T$ exactly $3-d_T(x)$ disjoint copies of $J$, each joined to $x$ by an edge at its subdivision vertex. Let $H=G_{K_1}$ (Figure~\ref{fig:H}).

\begin{proposition}\label{prop:GT}
For every tree $T$ with $p$ vertices and $\Delta(T)\le3$, the graph $G_T$ is connected and cubic, has $6p+10$ vertices, contains an induced copy of $J$, and has no perfect matching. Hence $\wdd(G_T)=4$.
\end{proposition}

\begin{proof}
Each vertex of $T$ receives $3-d_T(x)$ new edges, the subdivision vertex of each copy receives one, and the other vertices of a copy keep degree three. There are $\sum_x(3-d_T(x))=3p-2(p-1)=p+2$ copies, so $|V(G_T)|=p+5(p+2)=6p+10$. Each copy induces $J$. Removing $V(T)$ leaves $p+2$ odd components, more than $p$, so by Tutte's theorem there is no perfect matching, and Corollary~\ref{cor:cubic} gives $\wdd(G_T)=4$.
\end{proof}

\begin{proposition}[\citealp{Wa81,CRW86}]\label{prop:wallis}
Every connected cubic graph without a perfect matching has at least sixteen vertices, and $H$ is the only one with sixteen vertices.
\end{proposition}

\begin{proof}
Let $G$ be such a graph; $|V(G)|$ is even. By Tutte's theorem there is $U$ with $q=o(G-U)>|U|$, and $U\ne\emptyset$ because $G$ is connected of even order. By parity $q\ge|U|+2$. Each odd component $C$ sends $3|C|-2|E(C)|$ edges to $U$, an odd number. Let $\ell$ be the number of odd components sending exactly one edge; the others send at least three, so $3|U|\ge\ell+3(q-\ell)$, and $\ell\ge\frac32(q-|U|)\ge3$. A component sending one edge has one vertex of degree two and all others of degree three; it is not a single vertex, and it does not have three vertices, since its degree sum would be $8$, whereas a simple graph on three vertices has degree sum at most $6$. So it has at least five vertices, and $|V(G)|\ge|U|+5\ell\ge16$. Equality forces $|U|=1$, $\ell=3$ and three components of order five, each attached to the vertex of $U$ by one edge at its vertex of degree two. By Lemma~\ref{lem:33332} each is a copy of $J$, so $G\cong H$.
\end{proof}

\begin{theorem}[computer-assisted, via Theorem~\ref{thm:deg5}]\label{thm:GT}
Let $n=6p+10$ be the number of vertices of $G_T$. Every induced embedding of $G_T$ has degree at least six, and
\[
\eta(G_T)\ \ge\ \frac{3n}{2}.
\]
In particular $\eta(H)\ge24$, whereas the Erd\H{o}s--Kelly conditions with $d\ge4$ give only $\lceil5n/4\rceil$, which is $20$ for $H$. For $H$ itself, the bound $\eta(H)\ge24$ also follows without the computer (proof of Theorem~\ref{thm:exact}).
\end{theorem}

\begin{proof}
By Proposition~\ref{prop:GT}, $G_T$ is connected, contains an induced copy of $J$ and has $n\ge16$ vertices, so Corollary~\ref{cor:deg6} gives $d\ge6$. With $2m=3n$, condition \eqref{eq:c3} of Proposition~\ref{prop:EK} gives $N\ge2n-3n/d\ge2n-n/2=3n/2$. With $d\ge4$ alone, which is what Theorem~\ref{thm:nonab} gives, the same computation gives $5n/4$, and $N=\lceil5n/4\rceil$, $d=4$ satisfy \eqref{eq:c3}--\eqref{eq:c6}.
\end{proof}

\subsection{Explicit embeddings}

\begin{theorem}\label{thm:explicit}
Let
\[
\begin{aligned}
\Gamma_B&=\Z_2\times\Z_{10}, & S_B&=\{(0,\pm1),(1,\pm1),(1,\pm2)\},\\
\Gamma_H&=\Z_2\times\Z_{18}, & S_H&=\{(0,\pm4),(0,\pm7),(0,9),(1,0),(1,\pm4)\}.
\end{aligned}
\]
The labellings of Table~\ref{tab:labels} are induced embeddings of $B$ into $\Cay(\Gamma_B,S_B)$, of degree six, and of $H$ into $\Cay(\Gamma_H,S_H)$, of degree eight. Hence $\eta(B)\le20$ and $\eta(H)\le36$.
\end{theorem}

\begin{proof}
The labels are distinct. Computing the differences of the $16$ edges and $39$ non-edges of $B$, and of the $24$ edges and $96$ non-edges of $H$, in both orientations, gives $D_E=S_B$, respectively $D_E=S_H$, and $D_E\cap D_N=\emptyset$. Lemma~\ref{lem:diff} applies.
\end{proof}

\begin{table}[ht]
\centering
\small
\begin{tabular}{c ccccc c ccccc}
& \multicolumn{5}{c}{$B$ in $\Z_2\times\Z_{10}$; $\jv z\mapsto(1,4)$} && \multicolumn{5}{c}{$H$ in $\Z_2\times\Z_{18}$; $\jv z\mapsto(0,0)$}\\
$i$ & $\jv a_i$ & $\jv b_i$ & $\jv c_i$ & $\jv d_i$ & $\jv s_i$ && $\jv a_i$ & $\jv b_i$ & $\jv c_i$ & $\jv d_i$ & $\jv s_i$\\\midrule
1 & $(0,0)$ & $(1,2)$ & $(0,1)$ & $(1,1)$ & $(0,2)$ && $(1,2)$ & $(1,16)$ & $(0,2)$ & $(1,9)$ & $(0,9)$\\
2 & $(0,8)$ & $(1,6)$ & $(0,7)$ & $(1,7)$ & $(0,6)$ && $(1,1)$ & $(1,15)$ & $(0,15)$ & $(1,8)$ & $(0,4)$\\
3 &&&&&&& $(0,3)$ & $(0,17)$ & $(0,10)$ & $(1,3)$ & $(1,14)$
\end{tabular}
\caption{Induced embeddings of $B$ and $H$ (Theorem~\ref{thm:explicit}); the vertices of $J_i$ are named as in~\eqref{eq:EJ}.}
\label{tab:labels}
\end{table}

\begin{proposition}\label{prop:petersenhost}
Let $S=\{e_1,e_2,e_3,e_4,e_1+e_2+e_3+e_4\}\subseteq\F_2^4$. The ten elements of $\F_2^4$ that are neither $0$ nor adjacent to $0$ in $\Cay(\F_2^4,S)$ induce a Petersen graph. Hence $\eta(P)\le16$.
\end{proposition}

\begin{proof}
The ten elements are the vectors of weight two or three. Map a weight-two vector with support $\{i,j\}$ to $\{i,j\}$, and a weight-three vector missing $k\in\{1,2,3,4\}$ to $\{k,5\}$; this is a bijection onto the $2$-subsets of $\{1,\dots,5\}$. Two weight-two vectors are adjacent exactly when their supports are disjoint (their sum is $e_1+e_2+e_3+e_4$); a weight-two vector with support $\{i,j\}$ and a weight-three vector missing $k$ are adjacent exactly when $k\notin\{i,j\}$ (their sum has weight one); two weight-three vectors have sum of weight two and are not adjacent. So the induced graph is the Kneser graph $KG(5,2)$, the Petersen graph.
\end{proof}

\subsection{The minimum generator degree}

\begin{definition}\label{def:rab}
For a connected graph $G$ let $\rab(G)$ be the least degree $|S|$ of an induced embedding of $G$, over all orders.
\end{definition}

\begin{corollary}\label{cor:rab}
For every connected graph $G$, we have $\wdd(G)\le\rab(G)$ and $\eta(G)\ge\EK(G;\rab(G))$. The first inequality is strict for $B$ and for every $G_T$, which have $\wdd=4$ and $\rab\ge6$. Moreover, $\rab(B)=6$ and $6\le\rab(H)\le8$.
\end{corollary}

\begin{proof}
The first statement follows from Theorem~\ref{thm:nonab} and Definition~\ref{def:EK}. Strictness follows from Propositions~\ref{prop:B} and~\ref{prop:GT} and Theorem~\ref{thm:GT}, and the upper bounds from the hosts of degree six and eight in Theorem~\ref{thm:explicit}.
\end{proof}

\section{Exact values}\label{sec:exact}

\begin{theorem}\label{thm:exact}
For the Petersen graph $P$, the graph $B$, the Desargues graph $D$ and the cubic graph $H$,
\[
\eta(P)=16,\qquad\eta(B)=20,\qquad\eta(D)=32,\qquad\eta(H)=36.
\]
The abelian groups of order $\eta(G)$ that are hosts are $\Z_4\times\Z_4$, $\Z_2^2\times\Z_4$ and $\Z_2^4$ for $P$; $\Z_{20}$ and $\Z_2\times\Z_{10}$ for $B$; $\Z_2\times\Z_4\times\Z_4$, $\Z_2^3\times\Z_4$ and $\Z_2^5$ for $D$; and $\Z_2\times\Z_{18}$ for $H$.
\end{theorem}

The upper bounds follow from Propositions~\ref{prop:petersenhost} and~\ref{prop:desargues} and Theorem~\ref{thm:explicit}. We use a finite computation to obtain the group classification and the sharp lower bounds for $B$, $D$ and $H$. Section~\ref{sec:spectral} gives a separate proof of the sharp lower bound for $P$.

\subsection{Encoding and symmetry reductions}
We use the encoding of \citet[Section~8]{FB26c}. For a graph $G$ and a group $\Gamma$ there is a Boolean variable $x_{v,g}$ for each $v\in V$ and $g\in\Gamma$, meaning $\lambda(v)=g$, and a variable $y_c$ for each class $c=\{g,-g\}$ of non-zero elements, meaning $c\subseteq S$. The precise core CNF consists of the following clause families:
\begin{enumerate}
\item for each $v$, the at-least-one clause $\bigvee_{g\in\Gamma}x_{v,g}$ and the constraint that at most one of the $x_{v,g}$, $g\in\Gamma$, is true;
\item for each $g\in\Gamma$, the constraint that at most one of the $x_{v,g}$, $v\in V$, is true (injectivity);
\item for distinct $u,v$ and $g\ne h$, the clause $\neg x_{u,g}\vee\neg x_{v,h}\vee y_{c(g-h)}$ when $uv\in E$, and $\neg x_{u,g}\vee\neg x_{v,h}\vee\neg y_{c(g-h)}$ when $uv\notin E$.
\end{enumerate}
Here $c(g-h)=\{g-h,h-g\}$ is an inverse-pair class of nonzero group elements. Each at-most-one constraint on $k$ literals is written with the sequential counter of \citet{Si05}, which uses $k-1$ auxiliary variables and $3k-4$ clauses. The formula therefore has $n|\Gamma|$ label variables, one variable for each inverse-pair class, $n(|\Gamma|-1)+|\Gamma|(n-1)$ auxiliary variables and $O(n^2|\Gamma|^2)$ clauses. By Lemma~\ref{lem:diff} it is satisfiable exactly when $G$ has an induced embedding into a Cayley graph of $\Gamma$. We add the following symmetry-breaking clauses.

\begin{lemma}\label{lem:symmetry}
Let $(\Gamma,S,\lambda)$ be an induced embedding of $G$. For $g\in\Gamma$, $\varphi\in\operatorname{Aut}(\Gamma)$ and $\pi\in\operatorname{Aut}(G)$, the triples $(\Gamma,S,\lambda+g)$, $(\Gamma,\varphi(S),\varphi\circ\lambda)$ and $(\Gamma,S,\lambda\circ\pi)$ are induced embeddings. Fix a total order $\prec$ on $\Gamma$ and a set $R\subseteq\Gamma\setminus\{0\}$ meeting every orbit of $\operatorname{Aut}(\Gamma)$. If $G$ has an induced embedding into a Cayley graph of $\Gamma$, then it has one satisfying:
\begin{enumerate}
\item[(P)] $\lambda(p_0)=0$ and $\lambda(p_1)\in R$, for an edge $p_0p_1$ of the Petersen graph;
\item[(B)] for $B$ and for $B'$: $\lambda(\jv z)=0$, $\lambda(\jv s_1)\in R$, $\lambda(\jv a_i)\prec\lambda(\jv b_i)$ and $\lambda(\jv c_i)\prec\lambda(\jv d_i)$ for $i=1,2$;
\item[(H)] $\lambda(\jv z)=0$, $\lambda(\jv s_1)\in R$, $\lambda(\jv a_i)\prec\lambda(\jv b_i)$ and $\lambda(\jv c_i)\prec\lambda(\jv d_i)$ for $i=1,2,3$, and $\lambda(\jv s_2)\prec\lambda(\jv s_3)$; if $\Gamma=\Z_2^r$ with $r\ge3$, the conditions $\lambda(\jv s_1)\in R$ and $\lambda(\jv s_2)\prec\lambda(\jv s_3)$ may be replaced by $\lambda(\jv s_1)=e_1$, $\lambda(\jv s_2)=e_2$ and $\lambda(\jv s_3)\in\{e_1+e_2,e_3\}$;
\item[(D)] with the vertices of $D$ written as the $2$- and $3$-subsets of $\{1,\dots,5\}$: $\lambda(\{1,2\})=0$, $\lambda(\{1,2,3\})\in R$, $\lambda(\{1,2,4\})\prec\lambda(\{1,2,5\})$ and $\lambda(\{1,3\})\prec\lambda(\{2,3\})$.
\end{enumerate}
\end{lemma}

\begin{proof}
Translations preserve differences; $\varphi$ maps the differences of $\lambda$ to those of $\varphi\circ\lambda$ and $S$ to the connection set $\varphi(S)$; and $\pi$ preserves edges and non-edges. In each case, translate the first named vertex to $0$, and apply an automorphism of $\Gamma$, which fixes $0$, to move the label of the second into $R$. The remaining conditions are obtained with automorphisms of $G$ that fix these two vertices; hence these later automorphisms preserve the earlier translation and group-orbit normalisations. In $B$ and $H$, $\jv a_i,\jv b_i$ have the same closed neighbourhood and $\jv c_i,\jv d_i$ the same open neighbourhood, so the transpositions $(\jv a_i\,\jv b_i)$ and $(\jv c_i\,\jv d_i)$ are automorphisms; in $H$ the exchange of $J_2$ and $J_3$ is an automorphism fixing $\jv z$ and $J_1$, and it permutes the conditions on $J_2$ and $J_3$. If $\Gamma=\Z_2^r$, the labels of $\jv s_1,\jv s_2$ are distinct and non-zero, hence linearly independent over $\F_2$, and the label of $\jv s_3$ is either their sum or independent of both; a linear automorphism of $\Z_2^r$ gives the stated labels, and the transpositions are applied afterwards. In $D$, the transposition $(4\,5)$ of $\{1,\dots,5\}$ fixes $\{1,2\}$, $\{1,2,3\}$, $\{1,3\}$, $\{2,3\}$ and exchanges $\{1,2,4\}$ and $\{1,2,5\}$; then $(1\,2)$ fixes $\{1,2\}$, $\{1,2,3\}$, $\{1,2,4\}$, $\{1,2,5\}$ and exchanges $\{1,3\}$ and $\{2,3\}$.
\end{proof}

We take for $\prec$ the lexicographic order of coordinate vectors in invariant-factor form. For $R$ we take one element from each orbit of the subgroup of $\operatorname{Aut}(\Gamma)$ generated by multiplying a cyclic factor $\Z_{d_i}$ by a unit and by the maps $e_i\mapsto e_i+ce_j$ with $\operatorname{ord}(ce_j)\mid d_i$; every orbit of $\operatorname{Aut}(\Gamma)$ is a union of orbits of this subgroup, so $R$ meets all of them, whether or not the subgroup is all of $\operatorname{Aut}(\Gamma)$. Both kinds of map are automorphisms: multiplying $\Z_{d_i}$ by a unit is bijective, and $e_i\mapsto e_i+ce_j$ is well defined because $\operatorname{ord}(ce_j)\mid d_i$, with inverse $e_i\mapsto e_i-ce_j$. For every group used we also checked by computer that each of these maps is a bijective homomorphism.

\subsection{Certification}\label{sec:cert}
For each graph and each abelian group of the orders in Table~\ref{tab:computation}, we solved the formula with CaDiCaL 1.5.3 \citep{BFFH20} through PySAT \citep{IMM18}. Each satisfying assignment was turned into a labelling and checked with Lemma~\ref{lem:diff} on all pairs of vertices. For every unsatisfiable formula, CaDiCaL wrote a proof in binary DRAT format; it consists of lemmas with the reverse unit propagation property and of deletions, so it is a DRUP proof \citep{HHW13}. We checked every proof with the C version of the forward checker of \citet[Section~8]{FB26c}, which is independent of the solver and honours deletions. Every proof of at most one million lines was also checked with the Python version, $48$ of the $66$ proofs, and the two versions agreed on all of them. All $66$ proofs were verified; the longest, for $H$ in $\Z_{36}$, has $12.1$ million lines and took $23$ minutes to check. Table~\ref{tab:computation} summarises the computation. The labellings, checker outputs and SHA-256 digests of the formulas and proofs are in the deposit of \citet{FB26b}.

We regenerated all $76$ formulas and all $66$ unsatisfiability proofs. Every regenerated file matched its SHA-256 digest in \citet{FB26b}, and DRAT-trim \citep{WHH14} accepted all $66$ regenerated proofs. A separate implementation of the encoding reproduced the same clause multiset for each instance. It confirmed all $76$ outcomes, with DRAT-trim checking the negative cases and direct checks of the labellings for the ten positive cases.

\begin{table}[ht]
\centering
\small
\begin{tabular}{lcrcrrrrr}
\toprule
graph & orders & groups & order $\eta$: & \multicolumn{2}{c}{largest formula} & proof & check & solve\\
 & below $\eta$ & & non-hosts / hosts & vars & clauses & lines & (s) & (s)\\\midrule
Petersen $P$ & 14--15 & 2 & 2/3 & 469 & 11\,668 & 42\,k & 0.2 & 0.1\\
$B$ & 16--19 & 9 & 0/2 & 640 & 22\,949 & 115\,k & 0.4 & 0.4\\
Desargues $D$ & 28--31 & 5 & 4/3 & 1899 & 193\,190 & 14.4\,M & 743 & 274\\
$H$ & 24--35 & 24 & 3/1 & 1695 & 159\,128 & 48.2\,M & 3891 & 1017\\
$B'$ & 20--29 & 17 & 0/1 & 1142 & 71\,903 & 2.0\,M & 16 & 12\\\bottomrule
\end{tabular}
\caption{The computation behind Theorem~\ref{thm:exact} and Remark~\ref{rem:Bprime}. For each graph: the orders searched below $\eta$; the number of abelian groups of these orders, all unsatisfiable; the numbers of groups of order $\eta$ that are not and that are hosts; the size of the largest formula; the total length (lemmas and deletions) of the proofs of all unsatisfiable instances, each of which was verified; and the total times in seconds of the C checker (wall clock, on a shared two-core machine) and of the solver (processor time); both programs are single-threaded.}
\label{tab:computation}
\end{table}

\begin{proof}[Proof of Theorem~\ref{thm:exact}]
Corollary~\ref{cor:girth5} gives $\eta(P)\ge14$ and $\eta(D)\ge28$. For $B$ and $H$ the computation starts from $16$ and $24$, which need only Theorem~\ref{thm:nonab} and Corollary~\ref{cor:deg4}, not the computer-assisted Theorem~\ref{thm:deg5}. Indeed, Theorem~\ref{thm:nonab} gives $d\ge4$ for both, and Corollary~\ref{cor:deg4} excludes $d=4$. For $B$, $d\ge5$ gives $N\ge16$ as in the proof of Proposition~\ref{prop:B}. For $H$, a host of degree five has even order by Proposition~\ref{prop:EK}, so $N\ge2\cdot16-48/5$ forces $N\ge24$, and a host of degree at least six has $N\ge24$ by Corollary~\ref{cor:EKbound}. By Lemma~\ref{lem:symmetry}, the formulas of Table~\ref{tab:computation} for the orders below $\eta(G)$ are satisfiable whenever $G$ has an induced embedding of that order; they are all unsatisfiable, with verified proofs. The upper bounds were proved above, and the list of hosts comes from the formulas of order $\eta(G)$.
\end{proof}

For $P$ the computation is not needed: Theorem~\ref{thm:petersen} in Section~\ref{sec:spectral} proves $\eta(P)\ge16$ by hand, and Theorem~\ref{thm:desargues} proves $\eta(D)\ge30$.

\begin{remark}\label{rem:Bprime}
The same computation for the graph $B'$ of Example~\ref{ex:Bprime}, starting from its abelian bound $20$, shows that no abelian group of order $20$ to $29$ is a host (Table~\ref{tab:computation}). The labelling $\jv z\mapsto0$, $(\jv a_1,\jv b_1,\jv c_1,\jv d_1,\jv s_1)\mapsto(17,23,14,26,5)$, $(\jv a_2,\jv b_2,\jv c_2,\jv d_2,\jv s_2)\mapsto(7,13,4,16,25)$, $\ell_1\mapsto24$, $\ell_2\mapsto15$ is an induced embedding into $\Cay(\Z_{30},\{\pm3,\pm5,\pm6,\pm9\})$, of degree eight. Hence $\eta(B')=30$. Since $B$ is an induced subgraph of $B'$, in this example attaching a path of length two to $B$ raises $\eta$ from $20$ to $30$.
\end{remark}

\subsection{The limit of degree arguments}
The degree-based lower bounds of the preceding sections have a common form: a lower bound $d_0$ on the degree of every host, followed by the counting of Proposition~\ref{prop:EK}, which gives $\eta(G)\ge\EK(G;d_0)$. Such arguments cannot go beyond the least degree of a host.

\begin{proposition}\label{prop:limit}
If every host of $G$ has degree at least $d_0$, then
\[
\EK(G;d_0)\le\EK(G;\rab(G)).
\]
For the graphs of Table~\ref{tab:bounds},
\[
\begin{array}{lll}
P: & \rab(P)=5, & \EK(P;5)=14;\\
D: & \rab(D)=5, & \EK(D;5)=28;\\
B: & \rab(B)=6, & \EK(B;6)=17;\\
H: & 6\le\rab(H)\le8, & 24\le\EK(H;\rab(H))\le26;\\
B': & 6\le\rab(B')\le8, & 20\le\EK(B';\rab(B'))\le22.
\end{array}
\]
Hence the bounds $\eta(P)\ge14$, $\eta(D)\ge28$ and $\eta(B)\ge17$ are the best that any lower bound on the host degree can give, and for $H$ and $B'$ such bounds cannot exceed $26$ and $22$.
\end{proposition}

\begin{proof}
A host of degree $\rab(G)$ exists, so $d_0\le\rab(G)$, and $\EK(G;\cdot)$ is non-decreasing by Definition~\ref{def:EK}. The lower bounds on $\rab$ are Corollaries~\ref{cor:girth5} and~\ref{cor:deg6}, Proposition~\ref{prop:B}, Example~\ref{ex:Bprime} and Theorem~\ref{thm:GT}. The upper bounds come from hosts of degree five for $P$ and $D$ (Propositions~\ref{prop:petersenhost} and~\ref{prop:desargues}), six for $B$ and eight for $H$ (Theorem~\ref{thm:explicit}), and eight for $B'$ (Remark~\ref{rem:Bprime}). The values of $\EK$ follow from Definition~\ref{def:EK}. For example, for $H$ and $d=8$, condition \eqref{eq:c3} gives $t\ge(128-48)/8=10$, and $N=26$ satisfies \eqref{eq:c4}--\eqref{eq:c6}; for $d\ge9$, condition \eqref{eq:c3} gives $N\ge27$.
\end{proof}

Table~\ref{tab:bounds} collects the bounds. For $P$, $B$ and $D$ the degree method is exhausted, and the remaining gaps of $2$, $3$ and $4$ need another kind of argument. Section~\ref{sec:spectral} gives one, which closes the gap for $P$ and halves it for $D$. Condition (v) of \citet[Lemma~2.5]{FB26c}, which applies to the regular graphs, gives $13$, $25$ and $20$ for $P$, $D$ and $H$.

\begin{table}[ht]
\centering
\small
\begin{tabular}{lrrrrrrrr}
\toprule
graph & $n$ & $\EK(G;\Delta)$ & $\EK(G;\wdd)$ & $d_0$ & $\EK(G;d_0)$ & $\EK(G;\rab)$ & spectral & $\eta$\\\midrule
Petersen $P$ & 10 & 10 & 10 & 5 & 14 & 14 & 16 & 16\\
$B$ & 11 & 14 & 14 & 6 & 17 & 17 & -- & 20\\
$B'$ & 13 & 16 & 17 & 6 & 20 & 20--22 & -- & 30\\
Desargues $D$ & 20 & 20 & 20 & 5 & 28 & 28 & 30 & 32\\
$H$ & 16 & 16 & 20 & 6 & 24 & 24--26 & -- & 36\\
$G_T$ & $n$ & $n$ & $\lceil5n/4\rceil$ & 6 & $3n/2$ & & -- & unknown\\\bottomrule
\end{tabular}
\caption{Lower bounds and exact values. $\EK(G;\Delta)$ is the floor of \citet[Definition~2.6]{FB26c}; $\EK(G;\wdd)$ is Theorem~\ref{thm:nonab}; $d_0$ is the degree bound proved in Section~\ref{sec:abelian}, $\EK(G;\rab)$ is the limit of Proposition~\ref{prop:limit}, and the spectral column gives the bounds of Theorems~\ref{thm:petersen} and~\ref{thm:desargues}, which need no computer.}
\label{tab:bounds}
\end{table}

\begin{remark}\label{rem:brute}
Some structural results were also checked by brute force, over all connection sets of sizes three, four and five in the abelian groups of orders $5$ to $100$. No set of size three gives a Cayley graph containing an induced copy of $J$. Of the $822\,193$ sets of size four, $105$ do, all with $|\langle S\rangle|=7$; since every $4$-regular Cayley graph of $\Z_7$ is the complement of a $7$-cycle, $\Cay(\langle S\rangle,S)\cong C_7^2$ for all of them, as Theorem~\ref{thm:deg4} requires. Of the $8\,015\,695$ sets of size five, $44\,392$ do, all with $|\langle S\rangle|\in\{8,10,14\}$, as Theorem~\ref{thm:deg5} requires. The lattice membership test used for Theorem~\ref{thm:deg5} was compared, on $2000$ random integer lattices and $6000$ vectors, with an independent criterion based on the rank and the greatest common divisor of the maximal minors; the two agreed every time. The embeddings of Propositions~\ref{prop:desargues} and~\ref{prop:petersenhost}, Theorem~\ref{thm:explicit} and Remark~\ref{rem:Bprime} were checked with Lemma~\ref{lem:diff}.
\end{remark}

The SAT exclusions in Theorem~\ref{thm:exact} depend on the formulas of Table~\ref{tab:computation} and on Lemma~\ref{lem:symmetry}. The C and Python proof checkers are independent of CaDiCaL and have no formal verification. Section~\ref{sec:cert} also reports the DRAT-trim recheck against the exact regenerated formulas.

\section{Beyond degree arguments: a spectral obstruction}\label{sec:spectral}

By Proposition~\ref{prop:limit}, no lower bound on the host degree proves $\eta(P)=16$ or $\eta(D)>28$. At these values the counting of Proposition~\ref{prop:EK} is tight: every host vertex outside the image is forced into a rigid position. In an abelian Cayley graph this rigidity can be read off from the eigenvalues, which are character sums. This gives a computer-free proof of $\eta(P)=16$ and raises the computer-free bound for $D$ from $28$ to $30$.

Let $\widehat\Gamma$ be the group of characters $\chi\colon\Gamma\to\C^\times$, and write $\chi(S)=\sum_{s\in S}\chi(s)$; this number is real because $S=-S$. Let $A$ be the adjacency matrix of $X=\Cay(\Gamma,S)$, acting on functions $\Gamma\to\C$. Then $A\chi=\chi(S)\,\chi$ for every $\chi\in\widehat\Gamma$, and the characters form a basis \citep{Ba79}. We write $\ind_W$ for the indicator function of $W\subseteq\Gamma$ and $\ind=\ind_\Gamma$. A partition of the vertex set of a graph is \emph{equitable} if the number of neighbours that a vertex of one part has in another part depends only on the two parts; see \citet{GR01} and \citet{BH12} for the spectral properties of such partitions.

\begin{definition*}[{Tight embedding}]
Let $G$ be $k$-regular. An induced embedding $(\Gamma,S,\lambda)$ of $G$ of degree $d$ and order $N>n$ is \emph{tight} if equality holds in \eqref{eq:c3}, that is, $(N-n)d=n(d-k)$. A \emph{tight candidate pair} $(N,d)$ for $G$ is an integer pair with $N>n$ that satisfies \eqref{eq:c2} and \eqref{eq:c4}--\eqref{eq:c6}, with equality in \eqref{eq:c3}.
\end{definition*}

\begin{lemma}[Tight embeddings]\label{lem:tight}
Let $G$ be $k$-regular, let $(\Gamma,S,\lambda)$ be an induced embedding of $G$ of degree $d$ and order $N$, and put $t=N-n$ and $W=\Gamma\setminus\lambda(V)$. Suppose that $td=n(d-k)$, that is, equality holds in \eqref{eq:c3}, and $t\ge1$. Then $W$ is independent in $X$, and
\[
Af=(k-d)f\qquad\text{for}\qquad f=\ind_W-\tfrac tN\ind.
\]
Consequently some non-trivial character $\chi$ satisfies $\chi(S)=k-d$, and $\sum_{w\in W}\chi(w)=0$ for every non-trivial $\chi$ with $\chi(S)\ne k-d$.
\end{lemma}

\begin{proof}
Here $\sigma_d=nd-2m=n(d-k)$, so Proposition~\ref{prop:EK} gives $2|E(X[W])|=td-\sigma_d=0$. Hence every vertex of $W$ has all its $d$ neighbours in $U=\lambda(V)$, and every vertex of $U$ has exactly $d-k$ neighbours in $W$. Thus $A\ind_W=(d-k)(\ind-\ind_W)$. Since $A\ind=d\ind$ and $td/N=n(d-k)/N$,
\[
Af=(d-k)\ind-(d-k)\ind_W-\tfrac{n(d-k)}N\ind=(k-d)\Bigl(\ind_W-\tfrac{N-n}N\ind\Bigr)=(k-d)f.
\]
Expand $f=\sum_\chi c_\chi\chi$. Then $c_\chi(\chi(S)-k+d)=0$ for every $\chi$, the coefficient of the trivial character is the mean of $f$, which is $0$, and $c_\chi=N^{-1}\sum_{w\in W}\overline{\chi(w)}$ for $\chi\ne1$. Since $0<t<N$, $f\ne0$, so some $c_\chi\ne0$, and then $\chi\ne1$ and $\chi(S)=k-d$. Replacing $\chi$ by $\overline\chi$, which has the same value on $S$, gives the last statement.
\end{proof}

The second ingredient restricts which character sums can be integers.

\begin{lemma}[Rational character sums]\label{lem:galois}
Let $\chi\in\widehat\Gamma$ with $\chi(S)\in\Q$. For $s\in S$ let $m_s$ be the multiplicative order of $\chi(s)$. Then
\[
\chi(S)=\sum_{s\in S}\frac{\mu(m_s)}{\varphi(m_s)},
\]
where $\mu$ is the M\"obius function and $\varphi$ is Euler's function. In particular every term lies in $[-1,1]$, and equals $-1$ exactly when $\chi(s)=-1$.
\end{lemma}

\begin{proof}
Let $e$ be the exponent of $\Gamma$; the values of $\chi$ lie in $\Q(\zeta)$, $\zeta=e^{2\pi i/e}$. For $j$ prime to $e$ the field automorphism $\zeta\mapsto\zeta^j$ fixes $\chi(S)$, so $\chi(S)$ is the average over $j$ of $\sum_s\chi(s)^j$. If $\xi$ has order $m\mid e$, the reduction map $(\Z/e\Z)^\times\to(\Z/m\Z)^\times$ is onto with fibres of size $\varphi(e)/\varphi(m)$, so the average of $\xi^j$ is $\varphi(m)^{-1}$ times the sum of the primitive $m$-th roots of unity, which is $\mu(m)$ by M\"obius inversion of $\sum_{r\mid m}(\text{sum of primitive $r$-th roots})=[m=1]$. Finally $|\mu(m)/\varphi(m)|\le1$, with equality only for $m\in\{1,2\}$.
\end{proof}

In a group of odd order every term of Lemma~\ref{lem:galois} is at least $-\frac12$, with equality exactly when $m_s=3$. Together with Lemma~\ref{lem:tight} this classifies tight embeddings whose host degree is twice the degree of $G$. We need a fact about $\F_3^r$.

\begin{lemma}\label{lem:hyperplane}
Let $A\subseteq\F_3^r$ have $3^{r-1}$ elements, and let $v_1,\dots,v_r$ be a basis of $\F_3^r$ with $A\cap(A+v_i)=\emptyset$ for every $i$. Then $A$ is an affine hyperplane.
\end{lemma}

\begin{proof}
The linear automorphism of $\F_3^r$ that maps $v_i$ to the standard basis vector $e_i$ preserves the hypothesis and maps affine hyperplanes to affine hyperplanes, so we may take $v_i=e_i$. Any two points of a line $\{x,x+e_i,x+2e_i\}$ differ by $\pm e_i$, so each line meets $A$ at most once; there are $3^{r-1}=|A|$ lines in each direction, so each meets $A$ exactly once. For $r=1$, $A$ is a point. For $r\ge2$, $A$ is the graph of a map $F\colon\F_3^{r-1}\to\F_3$ that is a bijection in each argument when the others are fixed. Such a map is affine with non-zero coefficients. For one argument this holds because the six maps $x\mapsto ax+b$ with $a\ne0$ are all the permutations of $\F_3$. In general, write $F(x,z)=b(z)+\sum_ia_i(z)x_i$ by induction on the number of arguments. If $a_i(z)\ne a_i(z')$ for some $i$ and $z\ne z'$, the non-constant affine function $F(\cdot,z)-F(\cdot,z')$ has a zero, contradicting bijectivity in the last argument. So the $a_i$ are constant, and then $b$ is a permutation, hence affine.
\end{proof}

\begin{theorem}[Tight embeddings in groups of odd order]\label{thm:odd}
Let $G$ be $k$-regular and let $(\Gamma,S,\lambda)$ be an induced embedding of $G$ with $|\Gamma|$ odd, $\langle S\rangle=\Gamma$, $|S|=2k$ and equality in \eqref{eq:c3}. Then $W=\Gamma\setminus\lambda(V)$ is a coset of a subgroup $\Gamma_0$ of index three with $S\cap\Gamma_0=\emptyset$, and $G$ is isomorphic to the bipartite graph on $\Gamma_0\times\{0,1\}$ in which $(h,0)\sim(h',1)$ exactly when $h'-h\in Y$, for some $Y\subseteq\Gamma_0$ with $|Y|=k$. In particular $G$ is bipartite and $|\Gamma|=3n/2$.
\end{theorem}

\begin{proof}
With $d=2k$, equality in \eqref{eq:c3} reads $2kt=nk$, so $t=n/2$ and $|W|=N/3$. In particular $t\ge1$, so Lemma~\ref{lem:tight} applies, and $3$ divides $N$. By Lemma~\ref{lem:tight}, every non-trivial $\chi$ with $\sum_{w\in W}\chi(w)\ne0$ has $\chi(S)=-k$, a sum of $2k$ terms of Lemma~\ref{lem:galois}, each at least $-\frac12$. Hence $\chi(s)$ has order three for every $s\in S$, and $\chi^3=1$ because $S$ generates $\Gamma$. So $\ind_W$ is a combination of characters trivial on $3\Gamma$, and $W$ is a union of cosets of $3\Gamma$. Let $Q=\Gamma/3\Gamma\cong\F_3^r$ and let $A$ be the image of $W$, so $|A|=3^{r-1}$; since $|A|$ is an integer, $r\ge1$. The image of $S$ spans $Q$; choose $s_1,\dots,s_r\in S$ whose images $v_i$ form a basis. If $a$ and $a+v_i$ both lay in $A$, some $w\in W$ over $a$ would have $w+s_i\in W$, an edge inside the independent set $W$. By Lemma~\ref{lem:hyperplane}, $A$ is an affine hyperplane, so $W=a+\Gamma_0$ with $\Gamma_0$ of index three. No element of $S$ lies in $\Gamma_0$, since $W$ is independent. Thus $\lambda(V)$ is the union of the cosets $u_0+\Gamma_0$ and $u_1+\Gamma_0$, with $u_1-u_0=a_0\notin\Gamma_0$, and every edge joins them. Negation exchanges $a_0+\Gamma_0$ and $-a_0+\Gamma_0$, so $Y=\{s-a_0:s\in S\cap(a_0+\Gamma_0)\}$ has $k$ elements, and $(h,0)\mapsto u_0+h$, $(h,1)\mapsto u_1+h$ is the required isomorphism.
\end{proof}

Conversely, if $\Gamma_0$ has index three in an abelian group $\Gamma$ of odd order and $S=-S$ generates $\Gamma$ and misses $\Gamma_0$, then the graph induced on the complement of any coset of $\Gamma_0$ is a tight embedding of this kind. For example, $\Gamma=\Z_3\times\Z_7$, $\Gamma_0=\{0\}\times\Z_7$ and $S=\{\pm(1,y):y\in\{0,1,3\}\}$ give the Heawood graph.

In the sum of Lemma~\ref{lem:galois}, an inverse pair $\{a,-a\}$ contributes $2\mu(m_a)/\varphi(m_a)$, since $\chi(-a)=\overline{\chi(a)}$ has the same order, and an involution $h$ contributes $\chi(h)=\pm1$. For the divisors $m$ of $28$ the contribution of a pair is
\begin{equation}\label{eq:table28}
\begin{array}{c|cccccc}
m & 1 & 2 & 4 & 7 & 14 & 28\\\hline
2\mu(m)/\varphi(m) & 2 & -2 & 0 & -\frac13 & \frac13 & 0
\end{array}
\end{equation}

\begin{theorem}\label{thm:petersen}
$\eta(P)=16$ for the Petersen graph $P$, and no computer is needed for the lower bound.
\end{theorem}

\begin{proof}
Proposition~\ref{prop:petersenhost} gives $\eta(P)\le16$, and Corollary~\ref{cor:girth5} gives degree $d\ge5$ and order $N\ge14$. Suppose $N\in\{14,15\}$. Here $n=10$, $k=3$ and condition \eqref{eq:c3} reads $(N-10)d\ge10d-30$. For $N=14$ it gives $d\le5$, so $d=5$ and $4\cdot5=10\cdot2$. For $N=15$, $d$ is even by \eqref{eq:c6}, and \eqref{eq:c3} gives $d\le6$, so $d=6$ and $5\cdot6=10\cdot3$. In both cases Lemma~\ref{lem:tight} applies. A proper subgroup of $\Gamma$ has at most seven elements, so $\langle S\rangle=\Gamma$ by Lemma~\ref{lem:connected}.

\emph{Order $14$.} Then $\Gamma=\Z_{14}$, and since $|S|=5$ is odd, $S=\{\pm a,\pm b,h\}$ with $h=7$ the only involution. By Lemma~\ref{lem:tight} some $\chi$ has $\chi(S)=-2$. By Lemma~\ref{lem:galois} and \eqref{eq:table28}, the two pairs then contribute $-2-\chi(h)\in\{-3,-1\}$ in total, while each contributes one of $2,-2,-\frac13,\frac13$. No two of these values sum to $-3$ or $-1$.

\emph{Order $15$.} Then $|S|=6=2k$, the order is odd and equality holds in \eqref{eq:c3}. By Theorem~\ref{thm:odd}, the copy of $P$ would be bipartite. But $P$ contains five-cycles.
\end{proof}

\begin{theorem}\label{thm:desargues}
For the Desargues graph $D$, $\eta(D)\ge30$, without computer assistance.
\end{theorem}

\begin{proof}
By Corollaries~\ref{cor:C4free} and~\ref{cor:girth5}, $d\ge5$ and $N\ge28$. Here $n=20$, $k=3$ and \eqref{eq:c3} reads $(N-20)d\ge20d-60$. For $N=29$, $d$ is even, so $d\ge6$, but \eqref{eq:c3} gives $d\le60/11$; so $N\ne29$. For $N=28$ it gives $d\le5$, so $d=5$ with $8\cdot5=20\cdot2$, and Lemma~\ref{lem:tight} gives $\chi$ with $\chi(S)=-2$. The groups of order $28$ are $\Z_{28}$ and $\Z_2\times\Z_{14}$, of exponent dividing $28$, so \eqref{eq:table28} applies. Since $|S|$ is odd, $S$ contains one or three involutions. If $S=\{\pm a,\pm b,h\}$, the two pairs contribute $-2-\chi(h)\in\{-3,-1\}$, which no two of the values $2,-2,0,\pm\frac13$ sum to. If $S=\{\pm a,h_1,h_2,h_3\}$, the pair contributes $-2-\chi(h_1)-\chi(h_2)-\chi(h_3)\in\{-5,-3,-1,1\}$, which is none of these values.
\end{proof}

\begin{remark}\label{rem:scopeSpectral}
Lemma~\ref{lem:tight} needs a regular graph and equality in \eqref{eq:c3}, and Theorem~\ref{thm:odd} needs in addition odd order and host degree $2k$; we do not know a classification of tight embeddings in groups of even order. Throughout this remark we use the host degree bounds $d\ge5$ for $P$ and $D$ (Corollary~\ref{cor:C4free}) and $d\ge6$ for $H$ (Theorem~\ref{thm:GT}); without them, the definition would also admit the pairs $(25,4)$ for $D$ and $(20,4)$ for $H$. Below the respective values of $\eta$, the tight candidate pairs $(N,d)$ are $(14,5)$ and $(15,6)$ for $P$; $(28,5)$ and $(30,6)$ for $D$; and $(24,6)$, $(26,8)$, $(28,12)$ and $(29,16)$ for $H$; the graph $B$ is not regular. Theorems~\ref{thm:petersen} and~\ref{thm:desargues} use the pairs of $P$ and the pair $(28,5)$ of $D$. At $(30,6)$ Lemma~\ref{lem:tight} applies to $D$ and requires a character with $\chi(S)=-3$, but Theorem~\ref{thm:odd} does not, since the order is even; moreover $d=5$ also satisfies the counting conditions at order $30$ with strict inequality in \eqref{eq:c3}. At order $31$, $d=6$ is the only degree satisfying those conditions, again with strict inequality in \eqref{eq:c3}. For $H$ the four tight candidate pairs lie far below $\eta(H)=36$, and at the only odd order, $29$, the host degree is $16\ne2k$. Cauchy interlacing of the eigenvalues of $G$ with those of $X$ \citep{Ha95} gives further necessary conditions on $S$ without tightness, but where Lemma~\ref{lem:tight} applies it is weaker: for $D$ at order $28$, interlacing leaves $54$ of the $228$ generating connection sets of size five, and Lemma~\ref{lem:tight} leaves none.
\end{remark}

\section{Completions and induced embeddings}\label{sec:completions}

In the Cayley completion problem \citep{FB26a} the host has the same order as $G$. A \emph{completion certificate} for $G$ is a triple $(\Gamma,S,\lambda)$ with $|\Gamma|=n$, $\lambda\colon V\to\Gamma$ a bijection and $D_E(\lambda)\subseteq S$; its \emph{cost} is the number $n|S|/2-m$ of edges of $\Cay(\Gamma,S)$ that are not edges of $G$. A certificate of cost $0$ is an induced embedding of order $n$. Lemma~\ref{lem:classes} assumes only $D_E(\lambda)\subseteq S$, so a completion certificate also has a generator partition $(E_c)_{c\in\mathcal C(S)}$, and the pairs it pays for are
\[
P_c=\{uv\notin E:u\ne v,\ \lambda(u)-\lambda(v)\in c\},\qquad c\in\mathcal C(S).
\]
A completion certificate is an induced embedding exactly when every $P_c$ is empty. A secondary labelling can separate the real edges from the paid pairs. For $f\colon V\to\Z_{M}$ and a set $F$ of unordered pairs put $\delta_f(F)=\{\pm(f(u)-f(v)):uv\in F\}$.

\begin{theorem}[Transfer]\label{thm:transfer}
Let $(\Gamma,S,\lambda)$ be a completion certificate for $G$, let $M\ge1$, $f\colon V\to\Z_{M}$, and for each $c\in\mathcal C(S)$ let $J_c\subseteq\Z_{M}$ satisfy $J_c=-J_c$, $\delta_f(E_c)\subseteq J_c$ and $\delta_f(P_c)\cap J_c=\emptyset$. Then $\lambda'=(\lambda,f)\colon V\to\Gamma\times\Z_{M}$ and $S'=\{(s,j):s\in S,\ j\in J_{c(s)}\}$, where $c(s)$ is the class of $s$, form an induced embedding of $G$. In particular $\eta(G)\le nM$. Such sets $J_c$ exist exactly when $\delta_f(E_c)\cap\delta_f(P_c)=\emptyset$ for every $c$.
\end{theorem}

\begin{proof}
$\lambda'$ is injective because $\lambda$ is, and $S'=-S'$, $0\notin S'$ because $S$ and the $J_c$ are symmetric and $0\notin S$. Let $u\ne v$ and $s=\lambda(u)-\lambda(v)$. If $uv\in E$, then $s\in S$, $uv\in E_{c(s)}$ and $f(u)-f(v)\in J_{c(s)}$, so $\lambda'(u)-\lambda'(v)\in S'$. If $uv\notin E$ and $s\notin S$, no element of $S'$ has first coordinate $s$. If $uv\notin E$ and $s\in S$, then $uv\in P_{c(s)}$ and $f(u)-f(v)\notin J_{c(s)}$. The last sentence follows by taking $J_c=\delta_f(E_c)$.
\end{proof}

For $M=2$, the existence of $f$ reduces to the consistency of a linear system over $\F_2$, with one equation for each pair in $E_c\cup P_c$, restricted to the classes $c$ for which both $E_c$ and $P_c$ are non-empty. In the examples we know, Theorem~\ref{thm:transfer} does not improve a known value of $\eta$. The theorem gives an exact criterion for lifts with one symmetric second-coordinate set $J_c$ shared by all first-coordinate differences in each inverse class $c$. It uses the same generator partition as the lower bounds of Sections~\ref{sec:wd} and~\ref{sec:abelian}.

\section{Open problems}

\begin{problem}
Determine $\rab(H)$ and $\rab(B')$, which lie in $\{6,7,8\}$. More generally, determine $\eta(G_T)$ and $\rab(G_T)$ for the family of Proposition~\ref{prop:GT}.
\end{problem}

\begin{problem}
Proposition~\ref{prop:limit} shows that no lower bound on the host degree proves $\eta(P)=16$, $\eta(B)=20$ or $\eta(D)=32$. The spectral argument of Section~\ref{sec:spectral} proves the first and gives $\eta(D)\ge30$. Find a structural argument that proves $\eta(B)=20$ or $\eta(D)=32$, or that explains $\eta(H)=36$. In particular, can Lemma~\ref{lem:tight}, possibly combined with interlacing, exclude the tight candidate pair $(30,6)$ for $D$ (Remark~\ref{rem:scopeSpectral}), and which argument excludes the non-tight degree $d=5$ at order $30$?
\end{problem}

\begin{problem}
Call a graph $F$ \emph{nontrivially degree-forcing at level $d$} if $\Delta(F)\le d$ and there is $N$ such that no connected graph with at least $N$ vertices containing an induced copy of $F$ has an abelian host of degree $d$. By Theorems~\ref{thm:deg4} and~\ref{thm:deg5}, $J$ is degree-forcing at levels $3$, $4$ and $5$ but not at level $6$. Which graphs are nontrivially degree-forcing at some level above $\Delta(F)+2$? Is there, for every $d\ge3$, a graph of maximum degree at most $d$ that is nontrivially degree-forcing at level $d$?
\end{problem}

\begin{problem}
Given a partition of $E(G)$ into matchings and parts of maximum degree two, with prescribed generator orders, characterise when it is the generator partition of an induced abelian embedding.
\end{problem}

\begin{problem}
Lemma~\ref{lem:tight} turns a tight embedding of a $k$-regular graph into an equitable partition $(W,\lambda(V))$ of $\Cay(\Gamma,S)$ with quotient matrix $\bigl(\begin{smallmatrix}0&d\\d-k&k\end{smallmatrix}\bigr)$. For which $k$-regular graphs $G$ does such a partition exist with $\lambda(V)$ inducing a copy of $G$? Equivalently, when is the Erd\H{o}s--Kelly floor with equality in \eqref{eq:c3} attained by an abelian Cayley host?
\end{problem}

\section*{Use of AI tools}
The authors used AI assistant Claude Opus 5.5 (Anthropic) during the preparation of this paper for the computations of Section~\ref{sec:exact}, with the programs written for \citet{FB26c}; and to edit the text. AI assistant ChatGPT (Model GPT-6 Astra) were also used in the final revision, to review the text and to recheck the certificates of Section~\ref{sec:cert}. The authors retain responsibility for the mathematical correctness and for final validation of the statements, computations and reference metadata.

\section*{Funding}
This research received no specific grant from any funding agency in the public, commercial or not-for-profit sectors.

\section*{Competing interests}
The authors declare that they have no competing interests.

\section*{Data availability}
The data behind the computer-assisted results are in the deposit \citet{FB26b}, which is shared with \citet{FB26c}; its DOI is \href{https://doi.org/10.5281/zenodo.24457627}{10.5281/zenodo.24457627}. For Theorem~\ref{thm:deg5} and Remark~\ref{rem:brute} it contains the scripts \texttt{theorem45\_lattice.py}, \texttt{brute45\_groups.py} and \texttt{check\_hnf3.py} with their outputs; for Theorem~\ref{thm:sharp} the script \texttt{sharpness\_odd\_delta.py}; for the explicit embeddings of Section~\ref{sec:abelian} and Remark~\ref{rem:Bprime} the script \texttt{verify\_positive\_hosts.py}; and for Theorem~\ref{thm:exact} and Remark~\ref{rem:Bprime} the programs that write, solve and check the formulas, every labelling found, the SHA-256 digest of every formula and of every proof, and the output of the C checker on every proof and of the Python checker on every proof it checked. The deposit provides regeneration code for the formulas and proofs (about 0.8\,GB compressed). The authors will provide the files on request. The verification supplement of the recheck in Section~\ref{sec:cert}, with its scripts, the per-instance SHA-256 comparisons and the DRAT-trim outputs, is provided as ancillary files of the arXiv version of this paper.

\end{document}